\documentclass[a4paper,12pt]{article}

\usepackage{latexsym}
\usepackage{amssymb}
\usepackage{theorem}
\usepackage{amsmath}
\usepackage{amscd}
\usepackage{graphicx}
\usepackage{url}
\usepackage{color}

\newtheorem{theorem}{Theorem}[section]

\newtheorem{example}[theorem]{Example}
\newtheorem{proposition}[theorem]{Proposition}
\newtheorem{remark}[theorem]{Remark}
\newtheorem{definition}[theorem]{Definition}

\newcommand{\demo}{\par\noindent{\it Proof. \/}\ }
\newcommand{\enD}{\hfill $\Box$\vspace{3truemm} \par}
\newcommand{\R}{\mathbb{R}}
\newcommand{\N}{\mathbb{N}}

\newcommand{\bn}{\mbox{\boldmath $n$}}

\newcommand{\ba}{\mbox{\boldmath $a$}}
\newcommand{\bb}{\mbox{\boldmath $b$}}
\newcommand{\be}{\mbox{\boldmath $e$}}

\newcommand{\bv}{\mbox{\boldmath $v$}}
\newcommand{\bw}{\mbox{\boldmath $w$}}
\newcommand{\bx}{\mbox{\boldmath $x$}}
\newcommand{\by}{\mbox{\boldmath $y$}}
\newcommand{\bz}{\mbox{\boldmath $z$}}

\begin{document}

\title{Bertrand lightcone framed curves in the Lorentz-Minkowski 3-space}

\author{Nozomi Nakatsuyama and Masatomo Takahashi}

\date{\today}

\maketitle

\begin{abstract}
We consider mixed types of not only regular curves but also curves with singular points in the Lorentz-Minkowski 3-space.  
In order to consider mixed type of curves with singular points, we consider the lightcone frame and lightcone framed curves. 
By using lightcone frame, we can consider Bertrand types for lightcone framed curves, so-called Bertrand lightcone framed curves. 
We clarify that the existence conditions of Bertrand lightcone framed curves in all cases. 
As a consequence, we find pseudo-circular involutes and evolutes of mixed type curves which appear as Bertrand lightcone framed curves.
\end{abstract}

\renewcommand{\thefootnote}{\fnsymbol{footnote}}

\section*{Introduction}
Bertrand and Mannheim curves are classical objects in differential geometry (\cite{Aminov, Balgetir-Bektacs-Jun-ichi, Berger-Gostiaux, Bertrand, Chunxiao-Pei, HCIP, Papaioannou-Kiritsis, Struik}). 
A Bertrand (respectively, Mannheim) curve in the Euclidean $3$-space is a space curve whose principal normal line is the same as the principal normal (respectively, bi-normal) line of another curve. By definition, another curve is a parallel curve with respect to the direction of the principal normal vector. 
In \cite{Honda-Takahashi-2020}, they investigated the condition of the Bertrand and Mannheim curves of non-degenerate curves and framed curves. 
Moreover, we investigated the other cases, that is, a space curve whose tangent (or, principal normal, bi-normal) line is the same as the tangent (or, principal normal, bi-normal) line of another curve, respectively. 
We say that a Bertrand type curve if there exists such another curve. 
We investigated the existence conditions of Bertrand type curves in all cases in \cite{Nakatsuyama - Takahashi}. 
Moreover, we also investigated curves with singular points. As smooth curves with singular points, it is useful to use the framed curves in the Euclidean space (cf. \cite{Honda-Takahashi-2016}). 
We investigated the existence conditions of the Bertrand framed curves (Bertrand types of framed curves) in all cases in \cite{Nakatsuyama - Takahashi}. 
\par
For a non-lightlike non-degenerate curve, we have the arc-length parameter and the Frenet-Serret type formula by using a moving frame like a non-degenerate  space curve in the Euclidean space. 
It follows that we have the curvature of a non-lightlike non-degenerate curve (cf. \cite{López, O’Neil}).
\par
In \cite{Ucum-Ilarslan}, they investigated the necessary and sufficient conditions for timelike Bertrand curves in the Lorentz-Minkowski 3-space. 
In \cite{Hatice-Kazim}, they investigated the necessary and sufficient conditions for spacelike Bertrand curves with non-null normal vectors in the Lorentz-Minkowski 3-space.  
Also see in \cite{Wu-Zhou-Yao-Pei}. In \cite{Balgetir-Bektacs-Jun-ichi}, they investigated the necessary and sufficient conditions for null (lightlike) Bertrand curves in the Lorentz-Minkowski 3-space.  
\par
In this paper, we consider mixed types of not only regular curves but also curves with singular points in the Lorentz-Minkowski 3-space.  
In order to consider mixed types of curves with singular points, we consider the lightcone frame and lightcone framed curves (cf. \cite{Chen-Takahashi}).
By using the lightcone frame, we can consider Bertrand lightcone framed curves. 
We investigate a space curve whose lightlike (or, spacelike, timelike) vector is the same as the lightlike (or, spacelike, timelike) vector of another curve, respectively. 
We say that Bertrand lightcone framed curves  if there exists such another curve. 
In \S 3, we clarify that the existence conditions of Bertrand lightcone framed curves in all cases. 
As a consequence, we find psedo-circular evolutes and involutes of mixed type curves which appear as Bertrand lightcone framed curves (Theorems \ref{lightlike-mate-1}, \ref{lightlike-mate-2}, \ref{lightlike-mate-5}, \ref{lightlike-mate-7} and \ref{spacelike-mate4}). 
Therefore, it is useful to find new lightcone framed curves by using Bertrand lightcone framed curves. 
\par
We shall assume throughout the whole paper that all maps and manifolds are $C^{\infty}$ unless the contrary is explicitly stated.

\section{Preliminaries}
We review the theory of lightcone framed curves in the Lorentz-Minkowski 3-space. 
For details see in \cite{Chen-Takahashi}.

\subsection{Lorentz-Minkowski $3$-space}

The {\it Lorentz-Minkowski space} $\R^3_1$ is the space $\R^3$ endowed with the metric induced by the pseudo-scalar product $\langle \bx, \by \rangle=-x_1y_1+x_2y_2+x_3y_3$, where $\bx=(x_1,x_2,x_3)$ and $\by=(y_1,y_2,y_3)$. 
We say that a non-zero vector $\bx \in \R^3_1$ is {\it spacelike} if $\langle \bx, \bx \rangle >0$, {\it lightlike} if $\langle \bx, \bx \rangle=0$, and {\it timelike} if $\langle \bx, \bx \rangle <0$ respectively. 
The {\it norm} of a vector $\bx \in \R^3_1$ is defined by $||\bx||=\sqrt{|\langle \bx,\bx \rangle |}$. 
For any $\bx=(x_1, x_2, x_3), \by=(y_1, y_2, y_3)\in \R^3_1$, we define a vector $\bx\wedge \by$ by
$$
\bx\wedge \by= 
{\rm det} \left(\begin{array}{ccc}
-\be_1&\be_2&\be_3\\
x_1&x_2&x_3\\
y_1&y_2&y_3
\end{array}\right),
$$
where $\{\be_1, \be_2, \be_3\}$ is the canonical basis of $\R^3_1$. 
For any $\bw \in \R_1^3$, we can easily check that
$$
\langle \bw, \bx\wedge \by \rangle= \textrm{det}(\bw, \bx, \by),
$$
so that $\bx \wedge \by$ is pseudo-orthogonal to both $\bx$ and $\by$. 
Moreover, if $\bx$ is a unit timelike vector, $\by$ is a unit spacelike vector and $\langle \bx, \by\rangle=0$, $\bx\wedge \by=\bz$, then by a straightforward calculation we have
$\bz\wedge \bx=\by,\ \by\wedge \bz=-\bx.$
We define {\it hyperbolic $2$-space, de Sitter $2$-space} or {\it lightcone} by
\begin{align*}
H^2(-1) &=\{\bx\in \R_1^3\ |\ \langle\bx,\bx\rangle =-1\}, \\
S^2_1 &=\{\bx\in \R_1^3\ |\ \langle\bx,\bx\rangle =1\}, \\
LC^{*} &=\{\bx \in \R_1^3 \setminus \{0\}\ |\ \langle \bx, \bx\rangle=0\}.
\end{align*}
If $\bx=(x_0,x_1,x_2)$ is a non-zero lightlike vector, then $x_0 \not=0$.
Therefore, we have 
$$
\widetilde{\bx}=\left(1,\frac{x_1}{x_0},\frac{x_2}{x_0}\right) \in S^1_+=\{\bx=(x_0,x_1,x_2) \in LC^*| \ x_0=1\}.
$$
We call $S^1_+$ the {\it lightcone circle}.
In this paper, we consider two double Legendrian fibrations (cf. \cite{Izumiya09}):
\begin{itemize}
\item[(1)] 
$\Delta_1 =\{(\bv,\bw) \in H^2(-1) \times S^2_1 \ | \ \langle \bv,\bw \rangle=0 \},$
\par
$\pi_{11}:\Delta_1 \to H^2(-1), \pi_{11}(\bv,\bw)=\bv, \ \pi_{12}:\Delta_1 \to S^2_1, \pi_{12}(\bv,\bw)=\bw$, 
\par
$\theta_{11}=\langle d\bv, \bw \rangle|_{\Delta_1}, \theta_{12}=\langle \bv,d\bw \rangle|_{\Delta_1}$.
\par
\item[(2)] 
$\Delta_4 =\{(\bv,\bw) \in LC^* \times LC^* \ | \ \langle \bv,\bw \rangle=-2 \},$
\par
$\pi_{41}:\Delta_1 \to LC^*, \pi_{41}(\bv,\bw)=\bv, \ \pi_{42}:\Delta_1 \to LC^*, \pi_{42}(\bv,\bw)=\bw$, 
\par
$\theta_{41}=\langle d\bv, \bw \rangle|_{\Delta_4}, \theta_{42}=\langle \bv,d\bw \rangle|_{\Delta_4}$.
\end{itemize}
Note that $\Phi:\Delta_4 \to \Delta_1, \Phi(\bv,\bw)=((\bv+\bw)/2, (\bv-\bw)/2)$ is a contact diffeomorphsim. 
\par
Let $O(1,2)$ be the Lorentz group which consists of square matrices $A$ of order $3$ such that $ ^tAZA=Z$, where 
$$
Z:={\rm diag}(-1,1,1)=\left(\begin{array}{ccc}
-1 & 0 & 0 \\
0 & 1 & 0 \\
0 & 0 & 1
\end{array}\right).
$$
We set $SO(1,2)$ as $SO(1,2):=\{A \in O(1,2) \mid  {\rm det} A=1\}.$
For vectors $\ba, \bb \in \R^3_1$ and $A \in SO(1,2)$, we have 
$$
\langle \ba, \bb \rangle=\langle A(\ba), A(\bb) \rangle, \ 
A(\ba \wedge \bb)=A(\ba) \wedge A(\bb).
$$

\subsection{Lightcone framed curves}

In order to investigate mixed type of curves in the Lorentz-Minkowski space, we introduce the lightcone frame. 

\begin{definition}{\rm
Let $(\gamma,\ell^+,\ell^-):I \to \R^3_1 \times \Delta_4$ be a smooth mapping. 
We say that $(\gamma,\ell^+,\ell^-)$ is a {\it lightcone framed curve} if there exist smooth functions $\alpha, \beta:I \to \R$ such that $\dot{\gamma}(t)=\alpha(t) \ell^+(t)+\beta(t) \ell^-(t)$ for all $t \in I$. 
We also say that $\gamma$ is a {\it lightcone framed base curve} if there exists a smooth mapping $(\ell^+,\ell^-):I \to \Delta_4$ such that $(\gamma,\ell^+,\ell^-):I \to \R^3_1 \times \Delta_4$ is a lightcone  framed curve.
}
\end{definition}
Since $(\ell^+,\ell^-):I \to \Delta_4$, we have
$$
\langle \ell^+(t),\ell^+(t)\rangle=0, \ \langle \ell^-(t),\ell^-(t)\rangle=0, \ \langle \ell^+(t),\ell^-(t)\rangle=-2
$$
for all $t \in I$.
By $\dot{\gamma}(t)=\alpha(t) \ell^+(t)+\beta(t) \ell^-(t)$, $\langle \dot{\gamma}(t),\dot{\gamma}(t) \rangle=-4\alpha(t)\beta(t)$. 
Therefore, $\gamma$ is spacelike, lightlike, or timelike at $t$ if $\alpha(t)\beta(t)<0$, $\alpha(t)\beta(t)=0$ with $\alpha(t) \not=0$ or $\beta(t) \not=0$, or $\alpha(t)\beta(t)>0$, respectively.
Moreover, $t$ is a singular point of $\gamma$ if and only if $\alpha(t)=\beta(t)=0$.
\par
We denote $(\bn^T,\bn^S):I \to \Delta_1$,
$$
\bn^T(t)=\frac{\ell^+(t)+\ell^-(t)}{2}, \ \bn^S(t)=\frac{\ell^+(t)-\ell^-(t)}{2}.
$$
We define $\bn:I \to S^2_1, \bn(t)=\bn^T(t) \wedge \bn^S(t)=-(1/2)\ell^+(t) \wedge \ell^-(t)$. 
Then $\{\bn^T(t),\bn^S(t),\bn(t)\}$ is a pseudo-orthonormal frame of $\gamma(t)$.
We say that\\
$\{\ell^+(t),\ell^-(t),\bn(t)\}$ is a {\it lightcone frame} of $\gamma(t)$. 
Note that the lightcone frame is not a pseudo-orthonormal frame. 
By a direct calculation, we have 
\begin{eqnarray*}
\left(
\begin{array}{c}
\dot{\ell^{+}}(t)\\
\dot{\ell^{-}}(t)\\
\dot{\bn}(t)
\end{array}
\right)
&=&
\left(
\begin{array}{ccc}
\kappa_1(t)&0&2\kappa_3(t) \\
0&-\kappa_1(t)&2\kappa_2(t)\\
\kappa_2(t)&\kappa_3(t)&0
\end{array}
\right)
\left(
\begin{array}{c}
\ell^+(t)\\
\ell^-(t)\\
\bn(t)
\end{array}
\right), \\
\dot{\gamma}(t)&=&\alpha(t)\ell^+(t)+\beta(t)\ell^-(t). 
\end{eqnarray*}
We call $(\kappa_1,\kappa_2,\kappa_3,\alpha,\beta)$ a {\it (lightcone) curvature} of the lightcone framed curve $(\gamma,\ell^+,\ell^-):I \to \R^3_1 \times \Delta_4$.
\par
On the other hand, we have 
\begin{eqnarray*}\label{Frenet-type3}
\left(
\begin{array}{c}
\dot{\bn}^{T}(t)\\
\dot{\bn}^{S}(t)\\
\dot{\bn}(t)
\end{array}
\right)
&=&
\left(
\begin{array}{ccc}
0&\kappa_1(t)&\kappa^T(t) \\
\kappa_1(t)&0&-\kappa^S(t)\\
\kappa^T(t)&\kappa^S(t)&0
\end{array}
\right)
\left(
\begin{array}{c}
\bn^T(t)\\
\bn^S(t)\\
\bn(t)
\end{array}
\right), \\
\dot{\gamma}(t)&=&a(t)\bn^T(t)+b(t) \bn^S(t), \label{curve3}
\end{eqnarray*}
where $\kappa^T=\kappa_2+\kappa_3, \ \kappa^S=\kappa_2-\kappa_3, a=\alpha+\beta$ and $b=\alpha-\beta$.\par
In \cite{Chen-Takahashi}, we proved existence and uniqueness theorem for a special frame.
We can also prove the general case by the similar arguments, we omit it. 
\begin{theorem}[Existence theorem for lightcone framed curves]\label{existence}
Let $(\kappa_1,\kappa_2,\\\kappa_3,\alpha,\beta):$ $I \to \R^5$ be a smooth mapping. 
There exists a lightcone framed curve $(\gamma,\ell^+,\ell^-):I \to \R^3_1 \times \Delta_4$ such that the curvature is $(\kappa_1,\kappa_2,\kappa_3,\alpha,\beta)$.
\end{theorem}

\begin{definition}\label{congruence}{\rm
We say that $(\gamma,\ell^+,\ell^-)$ and $(\overline{\gamma},\overline{\ell}^+,\overline{\ell}^-):I \to \R^3_1 \times \Delta_4$ are {\it congruent as lightcone framed curves} if there exist $A \in SO(1,2)$ and $\ba \in \R^3_1$ such that $\overline{\gamma}(t)=A(\gamma(t))+\ba, \ \overline{\ell}^+(t)=A(\ell^+(t)), \ \overline{\ell}^-(t)=A(\ell^-(t))$ for all $t \in I$.
}
\end{definition}
\begin{theorem}[Uniqueness theorem for lightcone framed curves]\label{uniqueness}
Let $(\gamma,\ell^+,\\\ell^-)$ and $(\overline{\gamma},\overline{\ell}^+,\overline{\ell}^-):I \to \R^3_1 \times \Delta_4$ be lightcone framed curves with curvatures $(\kappa_1,\kappa_2,\kappa_3,\alpha,\beta)$ and $(\overline{\kappa}_1,\overline{\kappa}_2,\overline{\kappa}_3,\overline{\alpha},\overline{\beta})$, respectively.
$(\gamma,\ell^+,\ell^-)$ and $(\overline{\gamma},\overline{\ell}^+,\overline{\ell}^-)$ are congruent as lightcone framed curves if and only if $(\kappa_1,\kappa_2,\kappa_3,\alpha,\beta)=(\overline{\kappa}_1,\overline{\kappa}_2,\overline{\kappa}_3,\overline{\alpha},\overline{\beta})$.
\end{theorem}

\begin{proposition}\label{reflection}
Let $(\gamma,\ell^+,\ell^-):I \to \R^3_1 \times \Delta_4$ be a lightcone framed curve with curvature $(\kappa_1,\kappa_2,\kappa_3,\alpha,\beta)$.
Then $(\gamma,\ell^-,\ell^+)$ is also a lightcone framed curve with the curvature $(-\kappa_1,-\kappa_3,-\kappa_2,\beta,\alpha)$.
\end{proposition}

We consider a special lightcone frame.
We denote $\overline{\ell}^+(t)=(1/c(t))\ell^+(t)$ and  $\overline{\ell}^-(t)=c(t)\ell^-(t)$, where $c:I \to \R$ is a non-zero function. 
By a direct calculation, we have
$$
\overline{\ell}^+(t) \wedge \overline{\ell}^-(t)=\ell^+(t) \wedge \ell^-(t)=-2\bn(t). 
$$
Then $\{\overline{\ell}^+(t),\overline{\ell}^-(t),\bn(t) \}$ is also a lightcone frame of $\gamma(t)$.  
By a direct calculation, we have 
\begin{eqnarray*}
\left(
\begin{array}{c}
\dot{\overline{\ell}^{+}}(t)\\
\dot{\overline{\ell}^{-}}(t)\\
\dot{\bn}(t)
\end{array}
\right)
&=&
\left(
\begin{array}{ccc}
\overline{\kappa}_1(t)&0&2\overline{\kappa}_3(t) \\
0&-\overline{\kappa}_1(t)&2\overline{\kappa}_2(t)\\
\overline{\kappa}_2(t)&\overline{\kappa}_3(t)&0
\end{array}
\right)
\left(
\begin{array}{c}
\overline{\ell}^+(t)\\
\overline{\ell}^-(t)\\
\bn(t)
\end{array}
\right), \\
\dot{\gamma}(t)&=&\overline{\alpha}(t) \overline{\ell}^+(t)+\overline{\beta}(t) \overline{\ell}^-(t),
\end{eqnarray*}
where 
\begin{gather*}
\overline{\kappa}_1(t) =\frac{-\dot{c}(t)+c(t)\kappa_1(t)}{c^2(t)}, \ \overline{\kappa}_2(t)=c(t)\kappa_2(t), \ 
\overline{\kappa}_3(t)=\frac{\kappa_3(t)}{c(t)}, \\ 
\overline{\alpha}(t)=c(t)\alpha(t), \ 
\overline{\beta}(t)=\frac{\beta(t)}{c(t)}.
\end{gather*}
If we take $\dot{c}(t)=c(t)\kappa_1(t)$, that is, $c(t)=Ae^{\int \kappa_1(t) dt}$, where $A$ is a constant, then $\overline{\kappa}_1(t)=0$. 
Hence, we can always take $\overline{\kappa}_1(t)=0$.
We say that the lightcone frame $\{\overline{\ell}^+(t),\overline{\ell}^-(t),\bn(t) \}$ with $\overline{\kappa}_1(t)=0$ is an {\it adapted frame}. 

Moreover, we consider a special lightcone frame (cf. \cite{Liu-Pei}). 
Let $(\gamma,\ell^+,{\ell}^-):I \to \R^3_1 \times \Delta_4$ be a lightcone framed curve, where ${\ell}^+$ and $ {\ell}^-:I \to S^1_+$. 
Then there exists a smooth function $\theta:I \to \R$ such that 
$$
{\ell}^+(t)=(1,\cos \theta(t),\sin \theta(t)), \ 
{\ell}^-(t)=(1,-\cos \theta(t),-\sin \theta(t)).
$$
Therefore, 
$$
\bn(t)=-\frac{1}{2} {\ell}^+(t) \wedge {\ell}^-(t)=(0,-\sin \theta(t),\cos \theta(t)).
$$
We call the above lightcone frame $\{\ell^+(t),\ell^-(t),\bn(t)\}$ a {\it lightcone circle frame} of $\gamma(t)$. 
By a direct calculation, we have  
\begin{eqnarray*}
\left(
\begin{array}{c}
\dot{\ell^{+}}(t)\\
\dot{\ell^{-}}(t)\\
\dot{\bn}(t)
\end{array}
\right)
&=&
\left(
\begin{array}{ccc}
0&0&\dot{\theta}(t) \\
0&0&-\dot{\theta}(t)\\
-\dot{\theta}(t)/2&\dot{\theta}(t)/2&0
\end{array}
\right)
\left(
\begin{array}{c}
\ell^+(t)\\
\ell^-(t)\\
\bn(t)
\end{array}
\right), \\
\dot{\gamma}(t)&=&\alpha(t)\ell^+(t)+\beta(t)\ell^-(t).
\end{eqnarray*}
In this case, the curvature is given by $(0,-\dot{\theta}(t)/2,\dot{\theta}(t)/2,\alpha(t),\beta(t))$. 
In \cite{Liu-Pei}, $(\alpha,\beta,\theta)$ is called a {\it lightcone semi-polar coordinate}. 
and  evolutes of a mixed type curve are investigated by using lightcone circle frame of $\gamma$. 

\section{Bertrand types of lightcone framed curves}

Let $(\gamma,\ell^+,\ell^-)$ and $(\overline{\gamma},\overline{\ell}^+,\overline{\ell}^-):I \to \R^3_1 \times \Delta_4$ be lightcone framed curves with curvatures $(\kappa_1,\kappa_2,\kappa_3,\alpha,\beta)$ and $(\overline{\kappa}_1,\overline{\kappa}_2,\overline{\kappa}_3,\overline{\alpha},\overline{\beta})$. 
\begin{definition}{\rm
We say that $(\gamma,\ell^+,\ell^-)$ and $(\overline{\gamma},\overline{\ell}^+,\overline{\ell}^-)$ are {\it $(\bv,\overline{\bw})$-mates} if there exists a smooth function $\lambda:I \to \R$ with $\lambda \not\equiv 0$ such that $\overline{\gamma}(t)=\gamma(t)+\lambda(t)\bv(t)$ and $\bv(t)=\overline{\bw}(t)$ for all $t \in I$. 
We also say that $(\gamma,\ell^+,\ell^-)$ is a {\it $(\bv,\overline{\bw})$-Bertrand lightcone framed curve} (or, {\it $(\bv,\overline{\bw})$-Bertrand-Mannheim lightcone framed curve}) if there exists another lightcone framed curve $(\overline{\gamma},\overline{\ell}^+,\overline{\ell}^-)$ such that  $(\gamma,\ell^+,\ell^-)$ and $(\overline{\gamma},\overline{\ell}^+,\overline{\ell}^-)$ are $(\bv,\overline{\bw})$-mates. Here $\bv(t)$ and $\overline{\bw}(t)$ are non-zero vectors.
}
\end{definition}
Note that $\lambda \not\equiv 0$ means that $\{t \in I | \lambda(t) \not=0\}$ is a dense subset of $I$.
Then $\lambda$ is not identically zero for any non-trivial subintervals of $I$. 
It follows that $\gamma$ and $\overline{\gamma}$ are different space curves for any non-trivial subintervals of $I$.
By definition, $\bv$ and $\overline{\bw}$ are the same causal type. 
If the both vectors are lightlike, timelike or spacelike, then we say that $(\gamma,\ell^+,\ell^-)$ and $(\overline{\gamma},\overline{\ell}^+,\overline{\ell}^-)$ are {\it lightlike mate, timelike mate or spacelike mate}, respectively. 
Since $\ell^+$ and $\ell^-$ are lightlike, $\bn^T$ is a timelike, $\bn$ and $\bn^S$ are spacelike vectors, lightlike mates are four cases, timelike mate is one case, and spacelike mates are four cases. 
Moreover, we can consider naturally other lightlike vectors $\widetilde{\ell}^+=\bn^T+\bn$ and $\widetilde{\ell}^-=\bn^T-\bn$ by the moving frame $\{\bn^T,\bn^S,\bn\}$ of $\gamma$. 
We also add four cases in lightlike mates. 
By the construction, we have $(\widetilde{\ell}^+,\widetilde{\ell}^-):I \to \Delta_4$, however, $(\gamma,\widetilde{\ell}^+,\widetilde{\ell}^-)$ is not always a lightcone framed curve.
\par
We clarify that the existence conditions of Bertrand lightcone framed curves in all cases. 

\subsection{Lightlike mates}

Firstly, we consider lightlike mates. 
There are eight cases.
Let $(\gamma,\ell^+,\ell^-):I \to \R^3_1 \times \Delta_4$ be a lightcone framed curve with curvature $(\kappa_1,\kappa_2,\kappa_3,\alpha,\beta)$.

\begin{theorem}\label{lightlike-mate-1}
$(\gamma,\ell^+,\ell^-):I \to \R^3_1 \times \Delta_4$ is a $(\ell^+,\overline{\ell}^+)$-Bertrand lightcone framed curve if and only if there exist smooth functions $\lambda,k:I \to \R$ with $\lambda \not\equiv 0$ such that $k (t)\beta(t)=\lambda(t)\kappa_3(t)$ for all $t \in I$.
\end{theorem}
\demo
Suppose that $(\gamma,\ell^+,\ell^-):I \to \R^3_1 \times \Delta_4$ is a $(\ell^+,\overline{\ell}^+)$-Bertrand lightcone framed curve. Then there exists another lightcone framed curve $(\overline{\gamma},\overline{\ell}^+,\overline{\ell}^-)$ such that $(\gamma,\ell^+,\ell^-)$ and $(\overline{\gamma},\overline{\ell}^+,\overline{\ell}^-)$ are $(\ell^+,\overline{\ell}^+)$-mates, that is, 
there exists a smooth function $\lambda:I \to \R$ with $\lambda \not\equiv 0$ such that $\overline{\gamma}(t)=\gamma(t)+\lambda(t)\ell^+(t)$ and $\ell^+(t)=\overline{\ell}^+(t)$ for all $t \in I$. 
By differentiating $\overline{\gamma}(t)=\gamma(t)+\lambda(t)\ell^+(t)$, we have
$$
\overline{\alpha}(t)\overline{\ell}^+(t)+\overline{\beta}(t)\overline{\ell}^-(t)=(\alpha(t)+\dot{\lambda}(t)+\lambda(t)\kappa_1(t))\ell^+(t)+\beta(t)\ell^-(t)+2\lambda(t)\kappa_3(t)\bn(t).
$$
Since $\ell^+(t)=\overline{\ell}^+(t)$, we have $\overline{\beta}(t)=\beta(t)$. 
Moreover, there exist smooth functions $a_i,b_i,c_i:I \to \R, i=1,2$ such that 
\begin{align*}
\overline{\ell}^-(t) &=a_1(t)\ell^+(t)+b_1(t)\ell^{-}(t)+c_1(t)\bn(t),\\
\overline{\bn}(t) &=a_2(t)\ell^+(t)+b_2(t)\ell^{-}(t)+c_2(t)\bn(t),
\end{align*}
where $\overline{\bn}(t)=-(1/2) \overline{\ell}^+(t) \wedge \overline{\ell}^-(t)$. 
By the condition $(\overline{\ell}^+(t),\overline{\ell}^-(t)) \in \Delta_4$, we have 
$$
a_1(t)=a_2^2(t), \ b_1(t)=1, \ c_1(t)=2a_2(t), \ b_2(t)=0, \ c_2(t)=1.
$$
It follows that 
$$
\overline{\ell}^-(t) =a_2^2(t)\ell^+(t)+\ell^{-}(t)+2a_2(t)\bn(t), \ 
\overline{\bn}(t) =a_2(t)\ell^+(t)+\bn(t).
$$
Moreover, 
$$
\overline{\alpha}(t)\overline{\ell}^+(t)+\overline{\beta}(t)\overline{\ell}^-(t)=
(\overline{\alpha}(t)+\beta(t)a_2^2(t))\ell^+(t)+\beta(t)\ell^-(t)+2a_2(t)\bn(t).
$$
Therefore, we have $\overline{\alpha}(t)+\beta(t)a_2^2(t)=\alpha(t)+\dot{\lambda}(t)+\lambda(t)\kappa_1(t), \lambda(t)\kappa_3(t)=a_2(t)\beta(t)$. 
If we rewrite $a_2$ as $k$, then $k(t)\beta(t)=\lambda(t)\kappa_3(t)$ for all $t \in I$. 
\par
Conversely, suppose that there exist smooth functions $\lambda,k:I \to \R$ with $\lambda \not\equiv 0$ such that $k (t)\beta(t)=\lambda(t)\kappa_3(t)$ for all $t \in I$. 
If we consider $\overline{\gamma}(t)=\gamma(t)+\lambda(t)\ell^+(t), \overline{\ell}^{+}(t)=\ell^+(t),\overline{\ell}^-(t)=k^2(t)\ell^+(t)+\ell^-(t)+2k(t)\bn(t)$, 
\begin{align*}
\dot{\overline{\gamma}}(t) &=(\alpha(t)+\dot{\lambda}(t)+\lambda(t)\kappa_1(t))\ell^+(t)+\beta(t)\ell^-(t)+2\lambda(t)\kappa_3(t)\bn(t)\\
&=(\alpha(t)+\dot{\lambda}(t)+\lambda(t)\kappa_1(t)-\beta(t)k^2(t))\overline{\ell}^+(t)+\beta(t)\overline{\ell}^-(t)\\
&=\overline{\alpha}(t)\overline{\ell}^+(t)+\overline{\beta}(t)\overline{\ell}^-(t).
\end{align*}
It follows that $(\overline{\gamma},\overline{\ell}^+,\overline{\ell}^-):I \to \R^3_1 \times \Delta_4$ is a lightcone framed curve. 
Hence, $(\gamma,\ell^+,\ell^-)$ and $(\overline{\gamma},\overline{\ell}^+,\overline{\ell}^-)$ are $(\ell^+,\overline{\ell}^+)$-mates.
\enD

\begin{proposition}\label{curvature-lightlike-mate-1}
Suppose that there exist smooth functions $\lambda,k:I \to \R$ with $\lambda \not\equiv 0$ such that $k (t)\beta(t)=\lambda(t)\kappa_3(t)$ for all $t \in I$. 
Then the curvature $(\overline{\kappa}_1,\overline{\kappa}_2,\overline{\kappa}_3,\overline{\alpha},\overline{\beta})$ of the lightcone framed curve $(\overline{\gamma},\overline{\ell}^+,\overline{\ell}^-):I \to \R^3_1 \times \Delta_4$, where $\overline{\gamma}(t)=\gamma(t)+\lambda(t)\ell^+(t), \overline{\ell}^{+}(t)=\ell^+(t),\overline{\ell}^-(t)=k^2(t)\ell^+(t)+\ell^-(t)+2k(t)\bn(t)$ is given by 
\begin{gather*}
\overline{\kappa}_1(t)=\kappa_1(t)-2\kappa_3(t)k(t), \ 
\overline{\kappa}_2(t)=\dot{k}(t)+\kappa_1(t)k(t)+\kappa_2(t)-\kappa_3(t)k^2(t), \\
\overline{\kappa}_3(t)=\kappa_3(t), \
\overline{\alpha}(t)=\alpha(t)+\dot{\lambda}(t)+\lambda(t)\kappa_1(t)-\beta(t)k^2(t), \ 
\overline{\beta}(t)=\beta(t).
\end{gather*}
\end{proposition}
\demo
By a direct calculation, we have $\overline{\bn}(t)=k(t)\ell^+(t)+\bn(t)$ and 
\begin{align*}
\dot{\overline{\ell}^+}(t)&=\kappa_1(t)\ell^+(t)+2\kappa_3(t)\bn(t)\\
&=(\kappa_1(t)-2\kappa_3(t)k(t))\overline{\ell}^{+}+2\kappa_3(t)\overline{\bn}(t),\\
\dot{\overline{\ell}^-}(t)&=2k(t)\dot{k}(t)\ell^+(t)+k^2(t)(\kappa_1(t)\ell^+(t)+2\kappa_3(t)\bn(t))-\kappa_1(t)\ell^-(t)\\
&\quad +2\kappa_2(t)\bn(t)+2\dot{k}(t)\bn(t)+2k(t)(\kappa_2(t)\ell^+(t)+\kappa_3(t)\ell^{-}(t))\\
&=-(\kappa_1(t)-2\kappa_3(t)k(t))\overline{\ell}^-(t)\\
&\quad +2(\dot{k}(t)+\kappa_1(t)k(t)+\kappa_2(t)-\kappa_3(t)k^2(t))\overline{\bn}(t).
\end{align*}
Then we have 
$$
\overline{\kappa}_1(t)=\kappa_1(t)-2\kappa_3(t)k(t), \ 
\overline{\kappa}_2(t)=\dot{k}(t)+\kappa_1(t)k(t)+\kappa_2(t)-\kappa_3(t)k^2(t), \ 
\overline{\kappa}_3=\kappa_3(t).
$$
By the proof of Theorem \ref{lightlike-mate-1}, we also have 
$$
\overline{\alpha}(t)=\alpha(t)+\dot{\lambda}(t)+\lambda(t)\kappa_1(t)-\beta(t)k^2(t), \ 
\overline{\beta}(t)=\beta(t).
$$
\enD
\begin{remark}\label{lightlike-mate-1-re}{\rm 
$(1)$ Suppose that $\kappa_3(t) \not\equiv0$. If we take $k(t)=1$ and $\lambda(t)=\beta(t)/\kappa_3(t)$, then $(\gamma,\ell^+,\ell^-)$ is a $(\ell^+,\overline{\ell}^{+})$-Bertrand lightcone framed curve by Theorem \ref{lightlike-mate-1}.
In this case, $(\overline{\gamma},\overline{\ell}^+,\overline{\ell}^-)$ is given by 
$\overline{\gamma}(t)=\gamma(t)+(\beta(t)/\kappa_3(t))\ell^+(t),\\ \overline{\ell}^+(t)=\ell^+(t), \overline{\ell}^-(t)=\ell^+(t)+\ell^-(t)+2\bn(t)=2\widetilde{\ell}^+(t)$. 
It seems to be a kind of an evolute of mixed type curves with the direction $\ell^+$.\par
$(2)$ Suppose that $\beta(t) \not\equiv0$. If we take $k(t)=\kappa_2(t)$ and $\lambda(t)=\beta(t)$, then $(\gamma,\ell^+,\ell^-)$ is a $(\ell^+,\overline{\ell}^{+})$-Bertrand lightcone framed curve by Theorem \ref{lightlike-mate-1}.
In this case, $(\overline{\gamma},\overline{\ell}^+,\overline{\ell}^-)$ is given by 
$\overline{\gamma}(t)=\gamma(t)+\beta(t)\ell^+(t), \overline{\ell}^+(t)=\ell^+(t), \overline{\ell}^-(t)=\kappa_2^2(t)\ell^+(t)+\ell^-(t)+2\beta(t)\bn(t)$. 
}
\end{remark}

\begin{theorem}\label{lightlike-mate-2}
$(\gamma,\ell^+,\ell^-):I \to \R^3_1 \times \Delta_4$ is a $(\ell^+,\overline{\ell}^-)$-Bertrand lightcone framed curve if and only if there exist smooth functions $\lambda,k:I \to \R$ with $\lambda \not\equiv 0$ such that $k (t)\beta(t)=\lambda(t)\kappa_3(t)$ for all $t \in I$. 
That is, $(\gamma,\ell^+,\ell^-)$ is a $(\ell^+,\overline{\ell}^+)$-Bertrand lightcone framed curve if and only if $(\gamma,\ell^+,\ell^-)$ is a $(\ell^+,\overline{\ell}^-)$-Bertrand lightcone framed curve.
\end{theorem}
\demo
Suppose that $(\gamma,\ell^+,\ell^-):I \to \R^3_1 \times \Delta_4$ is a $(\ell^+,\overline{\ell}^-)$-Bertrand lightcone framed curve. Then there exists a lightcone framed curve $(\overline{\gamma},\overline{\ell}^+,\overline{\ell}^-):I \to \R^3_1 \times \Delta_4$ such that $(\gamma,\ell^+,\ell^-)$ and $(\overline{\gamma},\overline{\ell}^+,\overline{\ell}^-)$ are $(\ell^+,\overline{\ell}^-)$-mates.
By Proposition \ref{reflection}, $(\overline{\gamma},\overline{\ell}^-,\overline{\ell}^+)$ is also a lightcone framed curve. 
By Theorem \ref{lightlike-mate-1}, there exist smooth functions $\lambda,k:I \to \R$ with $\lambda \not\equiv 0$ such that $k (t)\beta(t)=\lambda(t)\kappa_3(t)$ for all $t \in I$.
The convese also holds by Theorem \ref{lightlike-mate-1} and Proposition \ref{reflection}.
\enD
By Propositions \ref{reflection} and \ref{curvature-lightlike-mate-1}, we have the following.
\begin{proposition}\label{curvature-lightlike-mate-2}
Suppose that there exist smooth functions $\lambda,k:I \to \R$ with $\lambda \not\equiv 0$ such that $k (t)\beta(t)=\lambda(t)\kappa_3(t)$ for all $t \in I$. 
Then the curvature $(\overline{\kappa}_1,\overline{\kappa}_2,\overline{\kappa}_3,\overline{\alpha},\overline{\beta})$ of the lightcone framed curve $(\overline{\gamma},\overline{\ell}^+,\overline{\ell}^-):I \to \R^3_1 \times \Delta_4$, where $\overline{\gamma}(t)=\gamma(t)+\lambda(t)\ell^+(t), \overline{\ell}^{+}(t)=k^2(t)\ell^+(t)+\ell^-(t)+2k(t)\bn(t),\overline{\ell}^-(t)=\ell^+(t)$ is given by 
\begin{align*}
&\overline{\kappa}_1(t)=-(\kappa_1(t)-2\kappa_3(t)k(t)), \ 
\overline{\kappa}_2(t)=-\kappa_3(t), \\ 
&\overline{\kappa}_3=-(\dot{k}(t)+\kappa_1(t)k(t)+\kappa_2(t)-\kappa_3(t)k^2(t)), \\
&\overline{\alpha}(t)=\beta(t), \ 
\overline{\beta}(t)=\alpha(t)+\dot{\lambda}(t)+\lambda(t)\kappa_1(t)-\beta(t)k^2(t).
\end{align*}
\end{proposition}
\begin{remark}{\rm
We can also prove Theorem \ref{lightlike-mate-2} and Proposition \ref{curvature-lightlike-mate-2} by the similar calculations of the proof of Theorem \ref{lightlike-mate-1} and Proposition \ref{curvature-lightlike-mate-1}. 
}
\end{remark}

\begin{theorem}\label{lightlike-mate-3}
$(\gamma,\ell^+,\ell^-):I \to \R^3_1 \times \Delta_4$ is a $(\ell^-,\overline{\ell}^+)$-Bertrand lightcone framed curve if and only if there exist smooth functions $\lambda,k:I \to \R$ with $\lambda \not\equiv 0$ such that $k (t)\alpha(t)=\lambda(t)\kappa_2(t)$ for all $t \in I$. 
\end{theorem}
\demo
Suppose that $(\gamma,\ell^+,\ell^-):I \to \R^3_1 \times \Delta_4$ is a $(\ell^-,\overline{\ell}^+)$-Bertrand lightcone framed curve. Then there exists another lightcone framed curve $(\overline{\gamma},\overline{\ell}^+,\overline{\ell}^-)$ such that $(\gamma,\ell^+,\ell^-)$ and $(\overline{\gamma},\overline{\ell}^+,\overline{\ell}^-)$ are $(\ell^-,\overline{\ell}^+)$-mates, that is, 
there exists a smooth function $\lambda:I \to \R$ with $\lambda \not\equiv 0$ such that $\overline{\gamma}(t)=\gamma(t)+\lambda(t)\ell^-(t)$ and $\ell^-(t)=\overline{\ell}^+(t)$ for all $t \in I$. 
By differentiating $\overline{\gamma}(t)=\gamma(t)+\lambda(t)\ell^-(t)$, we have
$$
\overline{\alpha}(t)\overline{\ell}^+(t)+\overline{\beta}(t)\overline{\ell}^-(t)=\alpha(t)\ell^+(t)+(\beta(t)+\dot\lambda(t)-\lambda(t)\kappa_1(t))\ell^-(t)+2\lambda(t)\kappa_2(t)\bn(t).
$$
Since $\ell^-(t)=\overline{\ell}^+(t)$, there exist smooth functions $a_i,b_i,c_i:I \to \R, i=1,2$ such that 
\begin{align*}
\overline{\ell}^-(t) &=a_1(t)\ell^+(t)+b_1(t)\ell^{-}(t)+c_1(t)\bn(t),\\
\overline{\bn}(t) &=a_2(t)\ell^+(t)+b_2(t)\ell^{-}(t)+c_2(t)\bn(t),
\end{align*}
where $\overline{\bn}(t)=-(1/2) \overline{\ell}^+(t) \wedge \overline{\ell}^-(t)$. 
By the condition $(\overline{\ell}^+(t),\overline{\ell}^-(t)) \in \Delta_4$, we have 
$$
a_1(t)=1, \ b_1(t)=b_2^2(t), \ c_1(t)=-2b_2(t), \ a_2(t)=0, \ c_2(t)=-1.
$$
It follows that 
$$
\overline{\ell}^-(t) =\ell^+(t)+b_2^2(t)\ell^{-}(t)-2b_2(t)\bn(t), \ 
\overline{\bn}(t) =b_2(t)\ell^-(t)-\bn(t).
$$
Moreover, 
$$
\overline{\alpha}(t)\overline{\ell}^+(t)+\overline{\beta}(t)\overline{\ell}^-(t)=
\overline{\beta}(t)\ell^+(t)+(\overline{\alpha}(t)+\overline{\beta}(t)b_2^2(t))\ell^-(t)-2b_2(t)\overline{\beta}(t)\bn(t).
$$
Therefore, we have $\overline{\beta}(t)=\alpha(t), \overline{\alpha}(t)+\overline{\beta}(t)b_2^2(t)=\beta(t)+\dot\lambda(t)-\lambda(t)\kappa_1(t),\\-b_2(t)\alpha(t)=\lambda(t)\kappa_2(t)$. 
If we rewrite $b_2$ as $-k$, then $k(t)\alpha(t)=\lambda(t)\kappa_2(t)$ for all $t \in I$. 
\par
Conversely, suppose that there exist smooth functions $\lambda,k:I \to \R$ with $\lambda \not\equiv 0$ such that $k(t)\alpha(t)=\lambda(t)\kappa_2(t)$ for all $t \in I$. 
If we consider $\overline{\gamma}(t)=\gamma(t)+\lambda(t)\ell^-(t), \overline{\ell}^{+}(t)=\ell^-(t),\overline{\ell}^-(t)=\ell^+(t)+k^2(t)\ell^-(t)+2k(t)\bn(t)$, 
\begin{align*}
\dot{\overline{\gamma}}(t) &=\alpha(t)\ell^+(t)+(\beta(t)+\dot\lambda(t)-\lambda(t)\kappa_1(t))\ell^-(t)+2\lambda(t)\kappa_2(t)\bn(t)\\
&=(\beta(t)+\dot\lambda(t)-\lambda(t)\kappa_1(t)-k^2(t)\alpha(t))\overline{\ell}^+(t)+\alpha(t)\overline{\ell}^-(t)\\
&=\overline{\alpha}(t)\overline{\ell}^+(t)+\overline{\beta}(t)\overline{\ell}^-(t).
\end{align*}
It follows that $(\overline{\gamma},\overline{\ell}^+,\overline{\ell}^-):I \to \R^3_1 \times \Delta_4$ is a lightcone framed curve. 
Hence, $(\gamma,\ell^+,\ell^-)$ and $(\overline{\gamma},\overline{\ell}^+,\overline{\ell}^-)$ are $(\ell^-,\overline{\ell}^+)$-mates.
\enD

\begin{proposition}\label{curvature-lightlike-mate-3}
Suppose that there exist smooth functions $\lambda,k:I \to \R$ with $\lambda \not\equiv 0$ such that $k (t)\alpha(t)=\lambda(t)\kappa_2(t)$ for all $t \in I$. 
Then the curvature $(\overline{\kappa}_1,\overline{\kappa}_2,\overline{\kappa}_3,\overline{\alpha},\overline{\beta})$ of the lightcone framed curve $(\overline{\gamma},\overline{\ell}^+,\overline{\ell}^-):I \to \R^3_1 \times \Delta_4$, where $\overline{\gamma}(t)=\gamma(t)+\lambda(t)\ell^-(t), \overline{\ell}^{+}(t)=\ell^-(t),\overline{\ell}^-(t)=\ell^+(t)+k^2(t)\ell^-(t)+2k(t)\bn(t)$ is given by 
\begin{align*}
&
\overline{\kappa}_1(t)=-\kappa_1(t)-2k(t)\kappa_2(t), \ 
\overline{\kappa}_2(t)=-\dot k(t)+k(t)\kappa_1(t)+k^2(t)\kappa_2(t)-\kappa_3(t), \\
&
\overline{\kappa}_3(t)=-\kappa_2(t), \ 
\overline{\alpha}(t)=\beta(t)+\dot\lambda(t)-\lambda(t)\kappa_1(t)-k^2(t)\alpha(t), \ 
\overline{\beta}(t)=\alpha(t).
\end{align*}
\end{proposition}
\demo
By a direct calculation, we have $\overline{\bn}(t)=-k(t)\ell^-(t)-\bn(t)$ and 
\begin{align*}
\dot{\overline{\ell}^+}(t)&=-\kappa_1(t)\ell^-(t)+2\kappa_2(t)\bn(t)\\
&=(-2k(t)\kappa_2(t)-\kappa_1(t))\overline{\ell}^{+}-2\kappa_2(t)\overline{\bn}(t),\\
\dot{\overline{\bn}}(t)&=-\kappa_2(t)\ell^+(t)+(-\dot k(t)+k(t)\kappa_1(t)-\kappa_3(t))\ell^-(t)-2k(t)\kappa_2(t)\bn\\
&=(-\dot k(t)+k(t)\kappa_1(t)+k^2(t)\kappa_2(t)-\kappa_3(t))\overline{\ell}^+(t)-\kappa_2(t)\overline{\ell}^-(t).
\end{align*}
Then we have 
\begin{align*}
&\overline{\kappa}_1(t)=-\kappa_1(t)-2k(t)\kappa_2(t), \\ 
&\overline{\kappa}_2(t)=-\dot k(t)+k(t)\kappa_1(t)+k^2(t)\kappa_2(t)-\kappa_3(t), \\ 
&\overline{\kappa}_3(t)=-\kappa_2(t).
\end{align*}
By the proof of Theorem \ref{lightlike-mate-3}, we also have 
$$
\overline{\alpha}(t)=\beta(t)+\dot\lambda(t)-\lambda(t)\kappa_1(t)-k^2(t)\alpha(t), \ 
\overline{\beta}(t)=\alpha(t).
$$
\enD
\begin{remark}\label{lightlike-mate-3-re}{\rm 
$(1)$ Suppose that $\kappa_2(t) \not\equiv0$. If we take $k(t)=1$ and $\lambda(t)=\alpha(t)/\kappa_2(t)$, then $(\gamma,\ell^+,\ell^-)$ is a $(\ell^-,\overline{\ell}^{+})$-Bertrand lightcone framed curve by Theorem \ref{lightlike-mate-3}.
In this case, $(\overline{\gamma},\overline{\ell}^+,\overline{\ell}^-)$ is given by 
$\overline{\gamma}(t)=\gamma(t)+(\alpha(t)/\kappa_2(t))\ell^-(t), \overline{\ell}^+(t)\\=\ell^-(t), \overline{\ell}^-(t)=\ell^+(t)+\ell^-(t)+2\bn(t)=2\widetilde{\ell}^+(t)$. 
It seems to be a kind of an evolute of mixed type curves with the direction $\ell^-$.\par
$(2)$ Suppose that $\alpha(t) \not\equiv0$. If we take $k(t)=\kappa_2(t)$ and $\lambda(t)=\alpha(t)$, then $(\gamma,\ell^+,\ell^-)$ is a $(\ell^-,\overline{\ell}^{+})$-Bertrand lightcone framed curve by Theorem \ref{lightlike-mate-3}.
In this case, $(\overline{\gamma},\overline{\ell}^+,\overline{\ell}^-)$ is given by 
$\overline{\gamma}(t)=\gamma(t)+\alpha(t)\ell^-(t), \overline{\ell}^+(t)=\ell^-(t), \overline{\ell}^-(t)=\ell^+(t)+\kappa_2^2(t)\ell^-(t)+2\kappa_2(t)\bn(t)$. 
}
\end{remark}

\begin{theorem}\label{lightlike-mate-4}
$(\gamma,\ell^+,\ell^-):I \to \R^3_1 \times \Delta_4$ is a $(\ell^-,\overline{\ell}^-)$-Bertrand lightcone framed curve if and only if there exist smooth functions $\lambda,k:I \to \R$ with $\lambda \not\equiv 0$ such that $k (t)\alpha(t)=\lambda(t)\kappa_2(t)$ for all $t \in I$. 
That is, $(\gamma,\ell^+,\ell^-)$ is a $(\ell^-,\overline{\ell}^+)$-Bertrand lightcone framed curve if and only if $(\gamma,\ell^+,\ell^-)$ is a $(\ell^-,\overline{\ell}^-)$-Bertrand lightcone framed curve.
\end{theorem}
\demo
Suppose that $(\gamma,\ell^+,\ell^-):I \to \R^3_1 \times \Delta_4$ is a $(\ell^-,\overline{\ell}^-)$-Bertrand lightcone framed curve. Then there exists a lightcone framed curve $(\overline{\gamma},\overline{\ell}^+,\overline{\ell}^-):I \to \R^3_1 \times \Delta_4$ such that $(\gamma,\ell^+,\ell^-)$ and $(\overline{\gamma},\overline{\ell}^+,\overline{\ell}^-)$ are  $(\ell^-,\overline{\ell}^+)$-mates.
By Proposition \ref{reflection}, $(\overline{\gamma},\overline{\ell}^-,\overline{\ell}^+)$ is also a lightcone framed curve. 
By Theorem \ref{lightlike-mate-3}, there exist smooth functions $\lambda,k:I \to \R$ with $\lambda \not\equiv 0$ such that $k (t)\alpha(t)=\lambda(t)\kappa_2(t)$ for all $t \in I$.
The convese also holds by Theorem \ref{lightlike-mate-3} and Proposition \ref{reflection}.
\enD

By Propositions \ref{reflection} and \ref{curvature-lightlike-mate-3}, we have the following.
\begin{proposition}\label{curvature-lightlike-mate-4}
Suppose that there exist smooth functions $\lambda,k:I \to \R$ with $\lambda \not\equiv 0$ such that $k (t)\alpha(t)=\lambda(t)\kappa_2(t)$ for all $t \in I$. 
Then the curvature $(\overline{\kappa}_1,\overline{\kappa}_2,\overline{\kappa}_3,\overline{\alpha},\overline{\beta})$ of the lightcone framed curve $(\overline{\gamma},\overline{\ell}^+,\overline{\ell}^-):I \to \R^3_1 \times \Delta_4$, where $\overline{\gamma}(t)=\gamma(t)+\lambda(t)\ell^-(t), \overline{\ell}^{+}(t)=\ell^+(t)+k^2(t)\ell^-(t)+2k(t)\bn(t),\overline{\ell}^-(t)=\ell^-(t)$ is given by 
\begin{align*}
&
\overline{\kappa}_1(t)=\kappa_1(t)+2k(t)\kappa_2(t), \ 
\overline{\kappa}_2(t)=\kappa_2(t), \\
&
\overline{\kappa}_3(t)=\dot k(t)-k(t)\kappa_1(t)-k^2(t)\kappa_2(t)+\kappa_3(t),\\
&
\overline{\alpha}(t)=\alpha(t), \ 
\overline{\beta}(t)=\beta(t)+\dot\lambda(t)-\lambda(t)\kappa_1(t)-k^2(t)\alpha(t).
\end{align*}
\end{proposition}
\begin{theorem}\label{lightlike-mate-5}
$(\gamma,\ell^+,\ell^-):I \to \R^3_1 \times \Delta_4$ is a $(\widetilde{\ell}^+, \overline{\ell}^+)$-Bertrand lightcone framed curve if and only if there exist smooth functions $\lambda, k:I \to \R$ with $\lambda \not\equiv 0$ such that $$k(t)(\alpha(t)+\beta(t))+\alpha(t)-\beta(t)+\lambda(t)(\kappa_1(t)+\kappa_2(t)-\kappa_3(t))=0$$ for all $t \in I$. 
\end{theorem}
\demo
Suppose that $(\gamma,\ell^+,\ell^-):I \to \R^3_1 \times \Delta_4$ is a $(\widetilde{\ell}^+, \overline{\ell}^+)$-Bertrand lightcone framed curve. Then there exists another lightcone framed curve $(\overline{\gamma},\overline{\ell}^+,\overline{\ell}^-)$ such that $(\gamma,\ell^+,\ell^-)$ and $(\overline{\gamma},\overline{\ell}^+,\overline{\ell}^-)$ are $(\widetilde{\ell}^+, \overline{\ell}^+)$-mates, that is, 
there exists a smooth function $\lambda:I \to \R$ with $\lambda \not\equiv 0$ such that $\overline{\gamma}(t)=\gamma(t)+\lambda(t)\widetilde{\ell}^+(t)$ and $\widetilde{\ell}^+(t)=\overline{\ell}^+(t)$ for all $t \in I$. 
By differentiating $\overline{\gamma}(t)=\gamma(t)+\lambda(t)\widetilde{\ell}^+(t)$, we have
\begin{align*}
\overline{a}(t)\overline{\bn}^T(t)+\overline{b}(t)\overline{\bn}^S(t)&=(\dot\lambda(t)+\lambda(t)\kappa^T(t))\bn(t)+(a(t)+\dot{\lambda}(t)+\lambda(t)\kappa^T(t))\bn^T(t)\\
&\quad+(b(t)+\lambda(t)(\kappa_1(t)+\kappa^S(t)))\bn^S(t).
\end{align*}
Since $\widetilde{\ell}^+(t)=\overline{\ell}^+(t)$, there exist smooth functions $a_i,b_i,c_i:I \to \R, i=1,2$ such that 
\begin{align*}
\overline{\ell}^-(t) &=a_1(t)\widetilde{\ell}^+(t)+b_1(t)\widetilde{\ell}^-(t)+c_1(t)\bn^S(t)\\
&=(a_1(t)+b_1(t))\bn^T(t)+c_1(t)\bn^S(t)+(a_1(t)-b_1(t))\bn(t),\\
\overline{\bn}(t) &=a_2(t)\widetilde{\ell}^+(t)+b_2(t)\widetilde{\ell}^-(t)(t)+c_2(t)\bn^S(t)\\
&=(a_2(t)+b_2(t))\bn^T(t)+c_2(t)\bn^S(t)+(a_2(t)-b_2(t))\bn(t).
\end{align*}
where $\overline{\bn}(t)=-(1/2) \overline{\ell}^+(t) \wedge \overline{\ell}^-(t)$. 
By the condition $(\overline{\ell}^+(t),\overline{\ell}^-(t)) \in \Delta_4$, we have 
$$
a_1(t)=a_2^2(t), \ b_1(t)=1, \ c_1(t)=-2a_2(t), \ b_2(t)=0, \ c_2(t)=-1.
$$
It follows that 
\begin{align*}
\overline{\ell}^-(t) &=(a_2^2(t)+1)\bn^T(t)-2a_2(t)\bn^S(t)+(a_2^2(t)-1)\bn(t), \\
&=a_2^2(t)\widetilde{\ell}^+(t)+\widetilde{\ell}^-(t)-2a_2(t)\bn^S(t),\\
\overline{\bn}(t) &=a_2(t)\bn^T(t)-\bn^S(t)+a_2(t)\bn(t) \\
&=a_2(t)\widetilde{\ell}^+(t)-\bn^S(t).
\end{align*}
Moreover, 
\begin{align*}
\overline{a}(t)\overline{\bn}^T(t)+\overline{b}(t)\overline{\bn}^S(t)&=(\overline{\alpha}(t)+(a_2^2(t)-1)\overline{\beta}(t))\bn(t)\\
&\quad+(\overline{\alpha}(t)+(a_2^2(t)+1)\overline{\beta}(t))\bn^T(t)
-2a_2(t)\overline{\beta}(t)\bn^S(t).  
\end{align*}
Therefore, we have $\overline{\alpha}(t)+(a_2^2(t)+1)\overline{\beta}(t)=\alpha(t)+\beta(t)+\dot{\lambda}(t)+\lambda(t)\kappa^T(t)$. 
If we rewrite $a_2$ as $k$, then $k(t)(\alpha(t)+\beta(t))+\alpha(t)-\beta(t)+\lambda(t)(\kappa_1(t)+\kappa_2(t)-\kappa_3(t))=0$ for all $t \in I$. 
\par
Conversely, suppose that there exist smooth functions $\lambda,k:I \to \R$ with $\lambda \not\equiv 0$ such that $k(t)(\alpha(t)+\beta(t))+\alpha(t)-\beta(t)+\lambda(t)(\kappa_1(t)+\kappa_2(t)-\kappa_3(t))=0$ for all $t \in I$. 
If we consider $\overline{\gamma}(t)=\gamma(t)+\lambda(t)\widetilde{\ell}^+(t), \overline{\ell}^{+}(t)=\widetilde{\ell}^+(t)=\bn^T(t)+\bn(t),\overline{\ell}^-(t)=(k^2(t)+1)\bn^T(t)-2k(t)\bn^S(t)+(k^2(t)-1)\bn(t)$, 
\begin{align*}
\dot{\overline{\gamma}}(t) &=(\dot\lambda(t)+\lambda(t)\kappa^T(t))\bn(t)+(\alpha(t)+\beta(t)+\dot{\lambda}(t)+\lambda(t)\kappa^T(t))\bn^T(t)\\
&\quad+(\alpha(t)-\beta(t)+\lambda(t)(\kappa_1(t)+\kappa^S(t)))\bn^S(t)\\
&=(a(t)+\dot{\lambda}(t)+\lambda(t)\kappa^T(t)-a(t)k^2(t)/2)\overline{\bn}^T(t)\\
&\quad+(\dot{\lambda}(t)+\lambda(t)\kappa^T(t)-a(t)k^2(t)/2)\overline{\bn}^S(t)\\
&=\overline{a}(t)\overline{\bn}^T(t)+\overline{b}(t)\overline{\bn}^S(t)\\
&=\frac{1}{2}(\overline{a}(t)+\overline{b}(t))\overline{\ell}^+(t)+\frac{1}{2}(\overline{a}(t)-\overline{b}(t))\overline{\ell}^-(t),
\end{align*}
where 
\begin{align*}
\overline{a}(t)=\dot\lambda(t)+\lambda(t)\kappa^T(t)-\frac{k^2(t)a(t)}{2}+a(t), \ \overline{b}(t)=\dot\lambda(t)+\lambda(t)\kappa^T(t)-\frac{k^2(t)a(t)}{2}.
\end{align*}
It follows that $(\overline{\gamma},\overline{\ell}^+,\overline{\ell}^-):I \to \R^3_1 \times \Delta_4$ is a lightcone framed curve. 
Hence, $(\gamma,\ell^+,\ell^-)$ and $(\overline{\gamma},\overline{\ell}^+,\overline{\ell}^-)$ are $(\widetilde{\ell}^+,\overline{\ell}^+)$-mates.
\enD
\begin{proposition}\label{curvature-lightlike-mate-5}
Suppose that there exist smooth functions $\lambda,k:I \to \R$ with $\lambda \not\equiv 0$ such that 
\begin{align*}
k(t)(\alpha(t)+\beta(t))+\alpha(t)-\beta(t)+\lambda(t)(\kappa_1(t)+\kappa_2(t)-\kappa_3(t))=0,
\end{align*}
for all $t \in I$. 
Then the curvature of the lightcone framed curve $(\overline{\gamma},\overline{\ell}^+,\overline{\ell}^-):I \to \R^3_1 \times \Delta_4$, where $\overline{\gamma}(t)=\gamma(t)+\lambda(t)\widetilde{\ell}^+(t), \overline{\ell}^{+}(t)=\widetilde{\ell}^+(t)=\bn^T(t)+\bn(t),\overline{\ell}^-(t)=(k^2(t)+1)\bn^T(t)-2k(t)\bn^S(t)+(k^2(t)-1)\bn(t)$ is given by 
$$
(\overline{\kappa}_1,\overline{\kappa}_2,\overline{\kappa}_3,\overline{\alpha},\overline{\beta})=(\overline{\kappa}_1,(1/2)(\overline{\kappa}^T+\overline{\kappa}^S),(1/2)(\overline{\kappa}^T-\overline{\kappa}^S),(1/2)(\overline{a}+\overline{b}),(1/2)(\overline{a}-\overline{b})),
$$
where
\begin{align*}
\overline{\kappa}_1(t)&=k(t)\kappa_1(t)+\kappa^T(t)+k(t)\kappa^S(t), \\
\overline{\kappa}^T(t)&=-\left(\frac{k^2(t)}{2}-1\right)\kappa_1(t)-k(t)\kappa^T(t)-\frac{k^2(t)}{2}\kappa^S(t)-\dot k(t), \\
\overline{\kappa}^S(t)&=\frac{k^2(t)}{2}\kappa_1(t)+k(t)\kappa^T(t)+\left(\frac{k^2(t)}{2}+1\right)\kappa^S(t)+\dot k(t),\\
\overline{a}(t)&=\dot\lambda(t)+\lambda(t)\kappa^T(t)-\frac{k^2(t)a(t)}{2}+a(t), \\ 
\overline{b}(t)&=\dot\lambda(t)+\lambda(t)\kappa^T(t)-\frac{k^2(t)a(t)}{2}.
\end{align*}
\end{proposition}
\demo
By a direct calculation,  
\begin{align*}
\overline{\bn}^T(t)&=\left(\frac{k^2(t)}{2}+1\right)\bn^T(t)-k(t)\bn^S(t)+\frac{k^2(t)}{2}\bn(t), \\
\overline{\bn}^S(t)&=-\frac{k^2(t)}{2}\bn^T(t)+k(t)\bn^S(t)+\left(-\frac{k^2(t)}{2}+1\right)\bn(t),\\
\overline{\bn}(t)&=k(t)\bn^T(t)-\bn^S(t)+k(t)\bn(t),\\
\dot{\overline{\bn}^T}(t)&=\left(\dot k(t)k(t)-k(t)\kappa_1(t)+\frac{k^2(t)}{2}\kappa^T(t)\right)\bn^T(t)\\
&\quad +\left(\left(\frac{k^2(t)}{2}+1\right)\kappa_1(t)-\dot k(t)+\frac{k^2(t)}{2}\kappa^S(t)\right)\bn^S(t)\\
&\quad +\left(\left(\frac{k^2(t)}{2}+1\right)\kappa^T(t)+k(t)\kappa^S(t)+\dot k(t)k(t)\right)\bn(t),\\
\dot{\overline{\bn}}(t)&=(\dot k(t)-\kappa_1(t)+k(t)\kappa^T(t))\bn^T(t)+k(t)(\kappa_1(t)+\kappa^S(t))\bn^S(t)\\
&\quad +(k(t)\kappa^T(t)+\kappa^S(t)+\dot k(t))\bn(t).
\end{align*}
Then we have 
\begin{align*}
\overline{\kappa}_1(t)&=k(t)\kappa_1(t)+\kappa^T(t)+k(t)\kappa^S(t), \\
\overline{\kappa}^T(t)&=-\left(\frac{k^2(t)}{2}-1\right)\kappa_1(t)-k(t)\kappa^T(t)-\frac{k^2(t)}{2}\kappa^S(t)-\dot k(t), \\
\overline{\kappa}^S(t)&=\frac{k^2(t)}{2}\kappa_1(t)+k(t)\kappa^T(t)+\left(\frac{k^2(t)}{2}+1\right)\kappa^S(t)+\dot k(t).
\end{align*}
By the proof of Theorem \ref{lightlike-mate-5}, we also have 
$$
\overline{a}(t)=\dot\lambda(t)+\lambda(t)\kappa^T(t)-\frac{k^2(t)a(t)}{2}+a(t), \ 
\overline{b}(t)=\dot\lambda(t)+\lambda(t)\kappa^T(t)-\frac{k^2(t)a(t)}{2}.
$$
\enD
\begin{remark}{\rm 
$(1)$ If $\kappa_1(t)+\kappa^S(t) \not=0, k(t)=0$ and $\lambda(t)=-b(t)/(\kappa_1(t)+\kappa^S(t))$, then $(\gamma,\ell^+,\ell^-)$ is a $(\widetilde{\ell}^+,\overline{\ell}^{+})$-Bertrand lightcone framed curve by Theorem \ref{lightlike-mate-5}.
In this case, $(\overline{\gamma},\overline{\ell}^+,\overline{\ell}^-)$ is given by 
$\overline{\gamma}(t)=\gamma(t)-(b(t)/(\kappa_1(t)+\kappa^S(t)))\widetilde{\ell}^+(t), \overline{\ell}^+(t)=\widetilde{\ell}^+(t), \overline{\ell}^-(t)=\bn^T(t)-\bn(t)=\widetilde{\ell}^-(t)$. 
\par
$(2)$ If $\kappa_1(t)+\kappa^S(t) \not=0, k(t)=1$ and $\lambda(t)=-2\alpha(t)/(\kappa_1(t)+\kappa^S(t))$, then $(\gamma,\ell^+,\ell^-)$ is a $(\widetilde{\ell}^+,\overline{\ell}^{+})$-Bertrand lightcone framed curve by Theorem \ref{lightlike-mate-5}.
In this case, $(\overline{\gamma},\overline{\ell}^+,\overline{\ell}^-)$ is given by 
$\overline{\gamma}(t)=\gamma(t)-(2\alpha(t)/(\kappa_1(t)+\kappa^S(t)))\widetilde{\ell}^+(t), \overline{\ell}^+(t)=\widetilde{\ell}^+(t), \overline{\ell}^-(t)=2\bn^T(t)-2\bn^S(t)=2{\ell}^-(t)$. 
\par
$(3)$ If $\kappa_1(t)+\kappa^S(t) \not=0, k(t)=-1$ and $\lambda(t)=2\beta(t)/(\kappa_1(t)+\kappa^S(t))$, then $(\gamma,\ell^+,\ell^-)$ is a $(\widetilde{\ell}^+,\overline{\ell}^{+})$-Bertrand lightcone framed curve by Theorem \ref{lightlike-mate-5}.
In this case, $(\overline{\gamma},\overline{\ell}^+,\overline{\ell}^-)$ is given by 
$\overline{\gamma}(t)=\gamma(t)+(2\beta(t)/(\kappa_1(t)+\kappa^S(t)))\widetilde{\ell}^+(t), \overline{\ell}^+(t)=\widetilde{\ell}^+(t), \overline{\ell}^-(t)=2\bn^T(t)+2\bn^S(t)=2{\ell}^+(t)$. 
\par
It seems to be a kind of psedo-circular evolutes of mixed type curves with the direction $\widetilde{\ell}^+$.
}
\end{remark}
\begin{theorem}\label{lightlike-mate-6}
$(\gamma,\ell^+,\ell^-):I \to \R^3_1 \times \Delta_4$ is a $(\widetilde{\ell}^+, \overline{\ell}^-)$-Bertrand lightcone framed curve if and only if there exist smooth functions $\lambda,k:I \to \R$ with $\lambda \not\equiv 0$ such that 
$$
k(t)(\alpha(t)+\beta(t))+\alpha(t)-\beta(t)+\lambda(t)(\kappa_1(t)+\kappa_2(t)-\kappa_3(t))=0
$$ 
for all $t \in I$. 
That is, $(\gamma,\ell^+,\ell^-)$ is a $(\widetilde{\ell}^+, \overline{\ell}^+)$-Bertrand lightcone framed curve if and only if $(\gamma,\ell^+,\ell^-)$ is a $(\widetilde{\ell}^+, \overline{\ell}^-)$-Bertrand lightcone framed curve.
\end{theorem}
\demo
Suppose that $(\gamma,\ell^+,\ell^-):I \to \R^3_1 \times \Delta_4$ is a $(\widetilde{\ell}^+, \overline{\ell}^-)$-Bertrand lightcone framed curve. Then there exists a lightcone framed curve $(\overline{\gamma},\overline{\ell}^+,\overline{\ell}^-):I \to \R^3_1 \times \Delta_4$ such that $(\gamma,\ell^+,\ell^-)$ and $(\overline{\gamma},\overline{\ell}^+,\overline{\ell}^-)$ are $(\widetilde{\ell}^+, \overline{\ell}^+)$-mates.
By Proposition \ref{reflection}, $(\overline{\gamma},\overline{\ell}^-,\overline{\ell}^+)$ is also a lightcone framed curve. 
By Theorem \ref{lightlike-mate-5}, there exist smooth functions $\lambda,k:I \to \R$ with $\lambda \not\equiv 0$ such that $k(t)(\alpha(t)+\beta(t))+\alpha(t)-\beta(t)+\lambda(t)(\kappa_1(t)+\kappa_2(t)-\kappa_3(t))=0$ for all $t \in I$.
The convese also holds by Theorem \ref{lightlike-mate-5} and Proposition \ref{reflection}.
\enD
By Propositions \ref{reflection} and \ref{curvature-lightlike-mate-5}, we have the following.
\begin{proposition}\label{curvature-lightlike-mate-6}
Suppose that there exist smooth functions $\lambda,k:I \to \R$ with $\lambda \not\equiv 0$ such that $$k(t)(\alpha(t)+\beta(t))+\alpha(t)-\beta(t)+\lambda(t)(\kappa_1(t)+\kappa_2(t)-\kappa_3(t))=0$$ for all $t \in I$. 
Then the curvature of the lightcone framed curve $(\overline{\gamma},\overline{\ell}^+,\overline{\ell}^-):I \to \R^3_1 \times \Delta_4$, where $\overline{\gamma}(t)=\gamma(t)+\lambda(t)\widetilde{\ell}^-(t), \overline{\ell}^{+}(t)=(k^2(t)+1)\bn^T(t)+2k(t)\bn^S(t)+(k^2(t)-1)\bn(t), \overline{\ell}^-(t)=\widetilde{\ell}^+(t)=\bn^T(t)+\bn(t)$ is given by 
$$
(\overline{\kappa}_1,\overline{\kappa}_2,\overline{\kappa}_3,\overline{\alpha},\overline{\beta})=(\overline{\kappa}_1,(1/2)(\overline{\kappa}^T+\overline{\kappa}^S),(1/2)(\overline{\kappa}^T-\overline{\kappa}^S),(1/2)(\overline{a}+\overline{b}),(1/2)(\overline{a}-\overline{b})),
$$
where
\begin{align*}
\overline{\kappa}_1(t)&=-k(t)\kappa_1(t)-\kappa^T(t)-k(t)\kappa^S(t), \\
\overline{\kappa}^T(t)&=\left(\frac{k^2(t)}{2}-1\right)\kappa_1(t)+k(t)\kappa^T(t)+\frac{k^2(t)}{2}\kappa^S(t)+\dot k(t), \\
\overline{\kappa}^S(t)&=\frac{k^2(t)}{2}\kappa_1(t)+k(t)\kappa^T(t)+\left(\frac{k^2(t)}{2}+1\right)\kappa^S(t)+\dot k(t),\\
\overline{a}(t)&=\dot\lambda(t)+\lambda(t)\kappa^T(t)-\frac{k^2(t)a(t)}{2}+a(t), \\ 
\overline{b}(t)&=-\dot\lambda(t)-\lambda(t)\kappa^T(t)+\frac{k^2(t)a(t)}{2}.
\end{align*}
\end{proposition}
\begin{theorem}\label{lightlike-mate-7}
$(\gamma,\ell^+,\ell^-):I \to \R^3_1 \times \Delta_4$ is a $(\widetilde{\ell}^-, \overline{\ell}^+)$-Bertrand lightcone framed curve if and only if there exist smooth functions $\lambda, k:I \to \R$ with $\lambda \not\equiv 0$ such that $$k(t)(\alpha(t)+\beta(t))-\alpha(t)+\beta(t)-\lambda(t)(\kappa_1(t)-\kappa_2(t)+\kappa_3(t))=0$$ for all $t \in I$. 
\end{theorem}
\demo
Suppose that $(\gamma,\ell^+,\ell^-):I \to \R^3_1 \times \Delta_4$ is a $(\widetilde{\ell}^-, \overline{\ell}^+)$-Bertrand lightcone framed curve. Then there exists another lightcone framed curve $(\overline{\gamma},\overline{\ell}^+,\overline{\ell}^-)$ such that $(\gamma,\ell^+,\ell^-)$ and $(\overline{\gamma},\overline{\ell}^+,\overline{\ell}^-)$ are $(\widetilde{\ell}^-, \overline{\ell}^+)$-mates, that is, 
there exists a smooth function $\lambda:I \to \R$ with $\lambda \not\equiv 0$ such that $\overline{\gamma}(t)=\gamma(t)+\lambda(t)\widetilde{\ell}^-(t)$ and $\widetilde{\ell}^-(t)=\overline{\ell}^+(t)$ for all $t \in I$. 
By differentiating $\overline{\gamma}(t)=\gamma(t)+\lambda(t)\widetilde{\ell}^-(t)$, we have
\begin{align*}
&\overline{a}(t)\overline{\bn}^T(t)+\overline{b}(t)\overline{\bn}^S(t)\\
&=(-\dot\lambda(t)+\lambda(t)\kappa^T(t))\bn(t)+(a(t)+\dot{\lambda}(t)-\lambda(t)\kappa^T(t))\bn^T(t)\\
&\quad+(b(t)+\lambda(t)(\kappa_1(t)-\kappa^S(t)))\bn^S(t).
\end{align*}
Since $\widetilde{\ell}^+(t)=\overline{\ell}^+(t)$, there exist smooth functions $a_i,b_i,c_i:I \to \R, i=1,2$ such that 
\begin{align*}
\overline{\ell}^-(t) &=a_1(t)\widetilde{\ell}^+(t)+b_1(t)\widetilde{\ell}^-(t)+c_1(t)\bn^S(t)\\
&=(a_1(t)+b_1(t))\bn^T(t)+c_1(t)\bn^S(t)+(a_1(t)-b_1(t))\bn(t),\\
\overline{\bn}(t) &=a_2(t)\widetilde{\ell}^+(t)+b_2(t)\widetilde{\ell}^-(t)(t)+c_2(t)\bn^S(t)\\
&=(a_2(t)+b_2(t))\bn^T(t)+c_2(t)\bn^S(t)+(a_2(t)-b_2(t))\bn(t).
\end{align*}
where $\overline{\bn}(t)=-(1/2) \overline{\ell}^+(t) \wedge \overline{\ell}^-(t)$. 
By the condition $(\overline{\ell}^+(t),\overline{\ell}^-(t)) \in \Delta_4$, we have 
$$
a_1(t)=1, \ a_2(t)=0, \ b_1(t)=b_2^2(t), \ c_1(t)=2b_2(t), \ c_2(t)=1.
$$
It follows that 
\begin{align*}
\overline{\ell}^-(t) &=(1+b_2^2(t))\bn^T(t)+2b_2(t)\bn^S(t)+(1-b_2^2(t))\bn(t), \\
&=\widetilde{\ell}^+(t)+b_2^2(t)\widetilde{\ell}^-(t)+2b_2(t)\bn^S(t),\\
\overline{\bn}(t) &=b_2(t)\bn^T(t)+\bn^S(t)-b_2(t)\bn(t) \\
&=b_2(t)\widetilde{\ell}^-(t)-\bn^S(t).
\end{align*}
Moreover, 
\begin{align*}
\overline{a}(t)\overline{\bn}^T(t)+\overline{b}(t)\overline{\bn}^S(t)&=(\overline{\alpha}(t)+(b_2^2(t)+1)\overline{\beta}(t))\bn(t)\\
&\quad-(\overline{\alpha}(t)+(b_2^2(t)+1)\overline{\beta}(t))\bn^T(t)+2b_2(t)\overline{\beta}(t)\bn^S(t).  
\end{align*}
Therefore, we have $\overline{\alpha}(t)+(b_2^2(t)+1)\overline{\beta}(t)=\alpha(t)+\beta(t)+\dot{\lambda}(t)-\lambda(t)\kappa^T(t)$. 
If we rewrite $b_2$ as $k$, then $k(t)(\alpha(t)+\beta(t))-\alpha(t)+\beta(t)-\lambda(t)(\kappa_1(t)-\kappa_2(t)+\kappa_3(t))=0$ for all $t \in I$. 
\par
Conversely, suppose that there exist smooth functions $\lambda,k:I \to \R$ with $\lambda \not\equiv 0$ such that $k(t)(\alpha(t)+\beta(t))-\alpha(t)+\beta(t)-\lambda(t)(\kappa_1(t)-\kappa_2(t)+\kappa_3(t))=0$ for all $t \in I$. 
If we consider $\overline{\gamma}(t)=\gamma(t)+\lambda(t)\widetilde{\ell}^-(t), \overline{\ell}^{+}(t)=\widetilde{\ell}^-(t)=\bn^T(t)-\bn(t),\overline{\ell}^-(t)=(1+k^2(t))\bn^T(t)+2k(t)\bn^S(t)+(1-k^2(t))\bn(t)$, 
\begin{align*}
\dot{\overline{\gamma}}(t) &=(-\dot\lambda(t)+\lambda(t)\kappa^T(t))\bn(t)+(\alpha(t)+\beta(t)+\dot{\lambda}(t)-\lambda(t)\kappa^T(t))\bn^T(t)\\
&\quad+(\alpha(t)-\beta(t)+\lambda(t)(\kappa_1(t)-\kappa^S(t)))\bn^S(t)\\
&=(a(t)+\dot{\lambda}(t)-\lambda(t)\kappa^T(t)-a(t)k^2(t)/2)\overline{\bn}^T(t)\\
&\quad+(\dot{\lambda}(t)-\lambda(t)\kappa^T(t)-a(t)k^2(t)/2)\overline{\bn}^S(t)\\
&=\overline{a}(t)\overline{\bn}^T(t)+\overline{b}(t)\overline{\bn}^S(t)\\
&=\frac{1}{2}(\overline{a}(t)+\overline{b}(t))\overline{\ell}^+(t)+\frac{1}{2}(\overline{a}(t)-\overline{b}(t))\overline{\ell}^-(t),
\end{align*}
where 
\begin{align*}
\overline{a}(t)=\dot\lambda(t)-\lambda(t)\kappa^T(t)-\frac{k^2(t)a(t)}{2}+a(t), \ \overline{b}(t)=\dot\lambda(t)-\lambda(t)\kappa^T(t)-\frac{k^2(t)a(t)}{2}.
\end{align*}
It follows that $(\overline{\gamma},\overline{\ell}^+,\overline{\ell}^-):I \to \R^3_1 \times \Delta_4$ is a lightcone framed curve. 
Hence, $(\gamma,\ell^+,\ell^-)$ and $(\overline{\gamma},\overline{\ell}^+,\overline{\ell}^-)$ are $(\widetilde{\ell}^-,\overline{\ell}^+)$-mates.
\enD

\begin{proposition}\label{curvature-lightlike-mate-7}
Suppose that there exist smooth functions $\lambda,k:I \to \R$ with $\lambda \not\equiv 0$ such that $$k(t)(\alpha(t)+\beta(t))-\alpha(t)+\beta(t)-\lambda(t)(\kappa_1(t)-\kappa_2(t)+\kappa_3(t))=0$$ for all $t \in I$. 
Then the curvature of the lightcone framed curve $(\overline{\gamma},\overline{\ell}^+,\overline{\ell}^-):I \to \R^3_1 \times \Delta_4$, where $\overline{\gamma}(t)=\gamma(t)+\lambda(t)\widetilde{\ell}^-(t), \overline{\ell}^{+}(t)=\widetilde{\ell}^-(t)=\bn^T(t)-\bn(t),\overline{\ell}^-(t)=(k^2(t)+1)\bn^T(t)+2k(t)\bn^S(t)+(-k^2(t)+1)\bn(t)$ is given by 
$$
(\overline{\kappa}_1,\overline{\kappa}_2,\overline{\kappa}_3,\overline{\alpha},\overline{\beta})=(\overline{\kappa}_1,(1/2)(\overline{\kappa}^T+\overline{\kappa}^S),(1/2)(\overline{\kappa}^T-\overline{\kappa}^S),(1/2)(\overline{a}+\overline{b}),(1/2)(\overline{a}-\overline{b})),
$$
where
\begin{align*}
\overline{\kappa}_1(t)&=-k(t)\kappa_1(t)-\kappa^T(t)+k(t)\kappa^S(t), \\
\overline{\kappa}^T(t)&=-\left(\frac{k^2(t)}{2}-1\right)\kappa_1(t)-k(t)\kappa^T(t)+\frac{k^2(t)}{2}\kappa^S(t)+\dot k(t), \\
\overline{\kappa}^S(t)&=-\frac{k^2(t)}{2}\kappa_1(t)-k(t)\kappa^T(t)+\left(\frac{k^2(t)}{2}+1\right)\kappa^S(t)+\dot k(t),\\
\overline{a}(t)&=\dot\lambda(t)-\lambda(t)\kappa^T(t)-\frac{k^2(t)a(t)}{2}+a(t), \\ 
\overline{b}(t)&=\dot\lambda(t)-\lambda(t)\kappa^T(t)-\frac{k^2(t)a(t)}{2}.
\end{align*}
\end{proposition}
\demo
By a direct calculation,  
\begin{align*}
\overline{\bn}^T(t)&=\left(\frac{k^2(t)}{2}+1\right)\bn^T(t)+k(t)\bn^S(t)-\frac{k^2(t)}{2}\bn(t), \\
\overline{\bn}^S(t)&=-\frac{k^2(t)}{2}\bn^T(t)-k(t)\bn^S(t)+\left(\frac{k^2(t)}{2}-1\right)\bn(t),\\
\overline{\bn}(t)&=k(t)\bn^T(t)+\bn^S(t)-k(t)\bn(t),\\
\dot{\overline{\bn}^T}(t)&=\left(\dot k(t)k(t)+k(t)\kappa_1(t)-\frac{k^2(t)}{2}\kappa^T(t)\right)\bn^T(t)\\
&\quad +\left(\left(\frac{k^2(t)}{2}+1\right)\kappa_1(t)+\dot k(t)-\frac{k^2(t)}{2}\kappa^S(t)\right)\bn^S(t)\\
&\quad +\left(\left(\frac{k^2(t)}{2}+1\right)\kappa^T(t)-k(t)\kappa^S(t)-\dot k(t)k(t)\right)\bn(t),\\
\dot{\overline{\bn}}(t)&=(\dot k(t)+\kappa_1(t)-k(t)\kappa^T(t))\bn^T(t)+k(t)(\kappa_1(t)-\kappa^S(t))\bn^S(t)\\
&\quad +(k(t)\kappa^T(t)-\kappa^S(t)-\dot k(t))\bn(t).
\end{align*}
Then we have 
\begin{align*}
\overline{\kappa}_1(t)&=-k(t)\kappa_1(t)-\kappa^T(t)+k(t)\kappa^S(t), \\
\overline{\kappa}^T(t)&=-\left(\frac{k^2(t)}{2}-1\right)\kappa_1(t)-k(t)\kappa^T(t)+\frac{k^2(t)}{2}\kappa^S(t)+\dot k(t), \\
\overline{\kappa}^S(t)&=-\frac{k^2(t)}{2}\kappa_1(t)-k(t)\kappa^T(t)+\left(\frac{k^2(t)}{2}+1\right)\kappa^S(t)+\dot k(t).
\end{align*}
By the proof of Theorem \ref{lightlike-mate-7}, we also have 
$$
\overline{a}(t)=\dot\lambda(t)-\lambda(t)\kappa^T(t)-\frac{k^2(t)a(t)}{2}+a(t), \ 
\overline{b}(t)=\dot\lambda(t)-\lambda(t)\kappa^T(t)-\frac{k^2(t)a(t)}{2}.
$$
\enD
\begin{remark}{\rm 
$(1)$ If $\kappa_1(t)-\kappa^S(t) \not=0, k(t)=0$ and $\lambda(t)=-b(t)/(\kappa_1(t)-\kappa^S(t))$, then $(\gamma,\ell^+,\ell^-)$ is a $(\widetilde{\ell}^-,\overline{\ell}^{+})$-Bertrand lightcone framed curve by Theorem \ref{lightlike-mate-7}.
In this case, $(\overline{\gamma},\overline{\ell}^+,\overline{\ell}^-)$ is given by 
$\overline{\gamma}(t)=\gamma(t)-(b(t)/(\kappa_1(t)-\kappa^S(t)))\widetilde{\ell}^-(t), \overline{\ell}^+(t)=\widetilde{\ell}^+(t), \overline{\ell}^-(t)=\bn^T(t)+\bn(t)=\widetilde{\ell}^+(t)$. 
\par
$(2)$ If $\kappa_1(t)-\kappa^S(t) \not=0, k(t)=1$ and $\lambda(t)=2\beta(t)/(\kappa_1(t)-\kappa^S(t))$, then $(\gamma,\ell^+,\ell^-)$ is a $(\widetilde{\ell}^-,\overline{\ell}^{+})$-Bertrand lightcone framed curve by Theorem \ref{lightlike-mate-7}.
In this case, $(\overline{\gamma},\overline{\ell}^+,\overline{\ell}^-)$ is given by 
$\overline{\gamma}(t)=\gamma(t)+(2\beta(t)/(\kappa_1(t)-\kappa^S(t)))\widetilde{\ell}^-(t), \overline{\ell}^+(t)=\widetilde{\ell}^-(t), \overline{\ell}^-(t)=2\bn^T(t)+2\bn^S(t)=2{\ell}^+(t)$. 
\par
$(3)$ If $\kappa_1(t)-\kappa^S(t) \not=0, k(t)=-1$ and $\lambda(t)=-2\alpha(t)/(\kappa_1(t)-\kappa^S(t))$, then $(\gamma,\ell^+,\ell^-)$ is a $(\widetilde{\ell}^-,\overline{\ell}^{+})$-Bertrand lightcone framed curve by Theorem \ref{lightlike-mate-7}.
In this case, $(\overline{\gamma},\overline{\ell}^+,\overline{\ell}^-)$ is given by 
$\overline{\gamma}(t)=\gamma(t)-(2\alpha(t)/(\kappa_1(t)-\kappa^S(t)))\widetilde{\ell}^-(t), \overline{\ell}^+(t)=\widetilde{\ell}^-(t), \overline{\ell}^-(t)=2\bn^T(t)-2\bn^S(t)=2{\ell}^-(t)$. 
\par
It seems to be a kind of psedo-circular evolutes of mixed type curves with the direction $\widetilde{\ell}^-$.
}
\end{remark}

\begin{theorem}\label{lightlike-mate-8}
$(\gamma,\ell^+,\ell^-):I \to \R^3_1 \times \Delta_4$ is a $(\widetilde{\ell}^-, \overline{\ell}^-)$-Bertrand lightcone framed curve if and only if there exist smooth functions $\lambda,k:I \to \R$ with $\lambda \not\equiv 0$ such that $$k(t)(\alpha(t)+\beta(t))-\alpha(t)+\beta(t)-\lambda(t)(\kappa_1(t)-\kappa_2(t)+\kappa_3(t))=0$$ for all $t \in I$. 
That is, $(\gamma,\ell^+,\ell^-)$ is a $(\widetilde{\ell}^-, \overline{\ell}^+)$-Bertrand lightcone framed curve if and only if $(\gamma,\ell^+,\ell^-)$ is a $(\widetilde{\ell}^-, \overline{\ell}^-)$-Bertrand lightcone framed curve.
\end{theorem}
\demo
Suppose that $(\gamma,\ell^+,\ell^-):I \to \R^3_1 \times \Delta_4$ is a $(\widetilde{\ell}^-, \overline{\ell}^-)$-Bertrand lightcone framed curve. Then there exists a lightcone framed curve $(\overline{\gamma},\overline{\ell}^+,\overline{\ell}^-):I \to \R^3_1 \times \Delta_4$ such that $(\gamma,\ell^+,\ell^-)$ and  $(\overline{\gamma},\overline{\ell}^+,\overline{\ell}^-)$ are $(\widetilde{\ell}^-, \overline{\ell}^-)$-mates.
By Proposition \ref{reflection}, $(\overline{\gamma},\overline{\ell}^-,\overline{\ell}^+)$ is also a lightcone framed curve. 
By Theorem \ref{lightlike-mate-7}, there exist smooth functions $\lambda,k:I \to \R$ with $\lambda \not\equiv 0$ such that $k(t)(\alpha(t)+\beta(t))-\alpha(t)+\beta(t)-\lambda(t)(\kappa_1(t)-\kappa_2(t)+\kappa_3(t))=0$ for all $t \in I$.
The convese also holds by Theorem \ref{lightlike-mate-7} and Proposition \ref{reflection}.
\enD
By Propositions \ref{reflection} and \ref{curvature-lightlike-mate-7}, we have the following.
\begin{proposition}\label{curvature-lightlike-mate-8}
Suppose that there exist smooth functions $\lambda,k:I \to \R$ with $\lambda \not\equiv 0$ such that $$k(t)(\alpha(t)+\beta(t))-\alpha(t)+\beta(t)-\lambda(t)(\kappa_1(t)-\kappa_2(t)+\kappa_3(t))=0$$ for all $t \in I$. 
Then the curvature of the lightcone framed curve $(\overline{\gamma},\overline{\ell}^+,\overline{\ell}^-):I \to \R^3_1 \times \Delta_4$, where $\overline{\gamma}(t)=\gamma(t)+\lambda(t)\widetilde{\ell}^-(t), \overline{\ell}^{+}(t)=(k^2(t)+1)\bn^T(t)-2k(t)\bn^S(t)-(k^2(t)-1)\bn(t), \overline{\ell}^-(t)=\widetilde{\ell}^-(t)=\bn^T(t)-\bn(t)$ is given by 
$$
(\overline{\kappa}_1,\overline{\kappa}_2,\overline{\kappa}_3,\overline{\alpha},\overline{\beta})=(\overline{\kappa}_1,(1/2)(\overline{\kappa}^T+\overline{\kappa}^S),(1/2)(\overline{\kappa}^T-\overline{\kappa}^S),(1/2)(\overline{a}+\overline{b}),(1/2)(\overline{a}-\overline{b})),
$$
where
\begin{align*}
\overline{\kappa}_1(t)&=k(t)\kappa_1(t)+\kappa^T(t)-k(t)\kappa^S(t), \\
\overline{\kappa}^T(t)&=\left(\frac{k^2(t)}{2}-1\right)\kappa_1(t)+k(t)\kappa^T(t)-\frac{k^2(t)}{2}\kappa^S(t)-\dot k(t), \\
\overline{\kappa}^S(t)&=-\frac{k^2(t)}{2}\kappa_1(t)-k(t)\kappa^T(t)+\left(\frac{k^2(t)}{2}+1\right)\kappa^S(t)+\dot k(t),\\
\overline{a}(t)&=\dot\lambda(t)-\lambda(t)\kappa^T(t)-\frac{k^2(t)a(t)}{2}+a(t), \\ 
\overline{b}(t)&=-\dot\lambda(t)+\lambda(t)\kappa^T(t)+\frac{k^2(t)a(t)}{2}.
\end{align*}
\end{proposition}

\subsection{Timelike mates}

Second, we consider only one timelike mate.

\begin{theorem}\label{timelike-mate}
$(\gamma,\ell^+,\ell^-):I \to \R^3_1 \times \Delta_4$ is a $(\bn^T,\overline{\bn}^T)$-Bertrand lightcone framed curve if and only if there exist smooth functions $\lambda,\theta:I \to \R$ with $\lambda \not\equiv 0$ such that 
$$
\lambda(t) (\kappa_2(t)+\kappa_3(t))\cos \theta(t)-(\alpha(t)-\beta(t)+\lambda(t)\kappa_1(t))\sin \theta(t)=0 
$$
for all $t \in I$. 
\end{theorem}
\demo
Suppose that $(\gamma,\ell^+,\ell^-):I \to \R^3_1 \times \Delta_4$ is a $(\bn^T,\overline{\bn}^T)$-Bertrand lightcone framed curve. Then there exists another lightcone framed curve $(\overline{\gamma},\overline{\ell}^+,\overline{\ell}^-)$ such that $(\gamma,\ell^+,\ell^-)$ and $(\overline{\gamma},\overline{\ell}^+,\overline{\ell}^-)$ are $(\bn^T,\overline{\bn}^T)$-mates, that is, 
there exists a smooth function $\lambda:I \to \R$ with $\lambda \not\equiv 0$ such that $\overline{\gamma}(t)=\gamma(t)+\lambda(t)\bn^T(t)$ and $\bn^T(t)=\overline{\bn}^T(t)$ for all $t \in I$. 
By differentiating $\overline{\gamma}(t)=\gamma(t)+\lambda(t) \bn^T(t)$, we have
$$
\overline{a}(t)\overline{\bn}^T(t)+\overline{b}(t)\overline{\bn}^S(t)=(a(t)+\dot{\lambda}(t))\bn^T(t)+(b(t)+\lambda(t)\kappa_1(t))\bn^S(t)+\lambda(t)\kappa^T(t)\bn(t).
$$
Since $\bn^T(t)=\overline{\bn}^T(t)$, we have $\overline{a}(t)=a(t)+\dot{\lambda}(t)$. 
Moreover, there exists a smooth function $\theta:I \to \R$ such that 
\begin{align*}
\begin{pmatrix}
\overline{\bn}(t)\\
\overline{\bn}^S(t)
\end{pmatrix}
=
\begin{pmatrix}
\cos \theta(t) & -\sin \theta(t)\\
\sin \theta(t) & \cos \theta(t)
\end{pmatrix}
\begin{pmatrix}
{\bn}(t)\\
{\bn}^S(t)
\end{pmatrix}.
\end{align*}
Then we have $\overline{b}(t)\sin \theta(t)=\lambda(t)\kappa^T(t)$ and 
$\overline{b}(t)\cos \theta(t)=b(t)+\lambda(t)\kappa_1(t)$.
It follows that $\overline{b}(t)=\lambda(t)\kappa^T(t) \sin \theta(t)+(b(t)+\lambda(t)\kappa_1(t))\cos \theta(t)$ and $\lambda(t)\kappa^T(t) \cos \theta(t)$\\
$-(b(t)+\lambda(t)\kappa_1(t))\sin \theta(t)=0$.
Therefore, we have
$$
\lambda(t) (\kappa_2(t)+\kappa_3(t))\cos \theta(t)-(\alpha(t)-\beta(t)+\lambda(t)\kappa_1(t))\sin \theta(t)=0
$$
for all $t \in I$.
\par
Conversely, suppose that there exist smooth functions $\lambda,\theta:I \to \R$ with $\lambda \not\equiv 0$ such that $\lambda(t)\kappa^T(t) \cos \theta(t)-(b(t)+\lambda(t)\kappa_1(t))\sin \theta(t)=0$ for all $t \in I$. 
If we consider $\overline{\gamma}(t)=\gamma(t)+\lambda(t)\bn^T(t), \overline{\bn}(t)=\cos \theta(t)\bn(t)-\sin \theta(t)\bn^S(t), \overline{\bn}^S(t)=\sin \theta(t)\bn(t)+\cos \theta(t)\bn^S(t)$, 
then $\overline{\bn}^T(t)=\bn^T(t)$ and  
\begin{align*}
\dot{\overline{\gamma}}(t) &=(a(t)+\dot{\lambda}(t))\bn^T(t)+(b(t)+\lambda(t)\kappa_1(t))\bn^S(t)+\lambda(t)\kappa^T(t)\bn(t)\\
&=(a(t)+\dot{\lambda}(t))\overline{\bn}^T(t)+(b(t)+\lambda(t)\kappa_1(t))(-\sin \theta(t) \overline{\bn}(t)+\cos \theta(t)\overline{\bn}^S(t))\\
&\quad +\lambda(t)\kappa^T(t)(\cos \theta(t) \overline{\bn}(t)+\sin \theta(t)\overline{\bn}^S(t)) \\
&=(a(t)+\dot{\lambda}(t))\overline{\bn}^T(t)+(\lambda(t)\kappa^T(t) \sin \theta(t)+(b(t)+\lambda(t)\kappa_1(t))\cos \theta(t))\overline{\bn}^S(t) \\
&= \overline{a}(t)\overline{\bn}^T(t)+\overline{b}(t)\overline{\bn}^S(t)\\
&=\frac{1}{2}(\overline{a}(t)+\overline{b}(t))\overline{\ell}^+(t)+\frac{1}{2}(\overline{a}(t)-\overline{b}(t))\overline{\ell}^-(t),
\end{align*}
where 
\begin{align*}
& \overline{a}(t)=a(t)+\dot{\lambda}(t), \ \overline{b}(t)=\lambda(t)\kappa^T(t) \sin \theta(t)+(b(t)+\lambda(t)\kappa_1(t))\cos \theta(t),\\
& \overline{\ell}^+(t)=\overline{\bn}^T(t)+\overline{\bn}^S(t), \ \overline{\ell}^-(t)=\overline{\bn}^T(t)-\overline{\bn}^S(t).
\end{align*}
It follows that $(\overline{\gamma},\overline{\ell}^+,\overline{\ell}^-):I \to \R^3_1 \times \Delta_4$ is a lightcone framed curve. 
Hence, $(\gamma,\ell^+,\ell^-)$ and $(\overline{\gamma},\overline{\ell}^+,\overline{\ell}^-)$ are $(\bn^T,\overline{\bn}^T)$-mates.
\enD

\begin{proposition}\label{curvature-timelike-mate}
Suppose that there exist smooth functions $\lambda,\theta:I \to \R$ with $\lambda \not\equiv 0$ such that
$$
\lambda(t) (\kappa_2(t)+\kappa_3(t))\cos \theta(t)-(\alpha(t)-\beta(t)+\lambda(t)\kappa_1(t))\sin \theta(t)=0 
$$ 
for all $t \in I$. 
Then the curvature of the lightcone framed curve $(\overline{\gamma},\overline{\ell}^+,\overline{\ell}^-):I \to \R^3_1 \times \Delta_4$, $\overline{\gamma}(t)=\gamma(t)+\lambda(t)\bn^T(t),\overline{\ell}^+(t)=\overline{\bn}^T(t)+\overline{\bn}^S(t), \overline{\ell}^-(t)=\overline{\bn}^T(t)-\overline{\bn}^S(t), 
\overline{\bn}^{T}(t)=\bn^T(t), 
\overline{\bn}^S(t)=\sin \theta(t)\bn(t)+\cos\theta(t) \bn^S(t)$ is given by 
$$
(\overline{\kappa}_1,\overline{\kappa}_2,\overline{\kappa}_3,\overline{\alpha},\overline{\beta})=(\overline{\kappa}_1,(1/2)(\overline{\kappa}^T+\overline{\kappa}^S),(1/2)(\overline{\kappa}^T-\overline{\kappa}^S),(1/2)(\overline{a}+\overline{b}),(1/2)(\overline{a}-\overline{b})),
$$
where
\begin{align*}
\overline{\kappa}_1(t)&=\kappa_1(t) \cos \theta(t)-\kappa^T(t)\sin \theta(t), \\ 
\overline{\kappa}^T(t)&=-\kappa_1(t)\sin \theta(t)+\kappa^T(t)\cos \theta(t), \\
\overline{\kappa}^S(t)&=\kappa^S(t)-\dot{\theta}(t),\\
\overline{a}(t)&=a(t)+\dot{\lambda}(t), \\ 
\overline{b}(t)&=\lambda(t)\kappa^T(t) \sin \theta(t)+(b(t)+\lambda(t)\kappa_1(t))\cos \theta(t).
\end{align*}
\end{proposition}
\demo
By a direct calculation,  
\begin{align*}
\dot{\overline{\bn}^T}(t)&=\kappa_1(t)\bn^S(t)+\kappa^T(t)\bn(t)\\
&=(\kappa_1(t) \cos \theta(t)+\kappa^T (t)\sin \theta(t))\overline{\bn}^{T}-(\kappa_1(t)\sin \theta(t)-\kappa^T(t)\cos \theta(t))\overline{\bn}(t), \\
\dot{\overline{\bn}^S}(t)&=\dot{\theta}(t)\cos \theta(t) \bn(t)+\sin \theta(t)(\kappa^T(t)\bn^T(t)+\kappa^S(t)\bn^S(t))\\
&\quad -\dot{\theta}(t)\sin \theta(t)\bn^S(t)+\cos \theta(t)(\kappa_1(t)\bn^T(t)-\kappa^S(t)\bn(t)) \\
&=(\kappa_1(t) \cos \theta(t)+\kappa^T(t) \sin \theta(t))\overline{\bn}^T(t)-(\kappa^S(t)-\dot{\theta}(t))\overline{\bn}(t).
\end{align*}
Then we have 
\begin{align*}
&\overline{\kappa}_1(t)=\kappa_1(t) \cos \theta(t)-\kappa^T(t)\sin \theta(t), \ 
\overline{\kappa}^T(t)=-\kappa_1(t)\sin \theta(t)+\kappa^T(t)\cos \theta(t), \\ 
&\overline{\kappa}^S(t)=\kappa^S(t)-\dot{\theta}(t).
\end{align*}
By the proof of Theorem \ref{timelike-mate}, we also have 
$$
\overline{a}(t)=a(t)+\dot{\lambda}(t), \ 
\overline{b}(t)=\lambda(t)\kappa^T(t) \sin \theta(t)+(b(t)+\lambda(t)\kappa_1(t))\cos \theta(t).
$$
\enD
\begin{remark}{\rm 
If $(\gamma,\ell^+,\ell^-):I \to \R^3_1 \times \Delta_4$ is a lightcone framed curve with lightcone circle frame, then $\kappa_1=0$ and $\kappa_2+\kappa_3=0$, see \S 2. 
By Theorem \ref{timelike-mate}, it is automatically holded, since we can take $\theta=0$. 
Hence $(\gamma,\ell^+,\ell^-)$ is a $(\bn^T,\overline{\bn}^T)$-Bertrand lightcone framed curve.
}
\end{remark}

\subsection{Spacelike mates}

Third, we consider spacelike mates. 
There are four cases.

\begin{proposition}\label{lambda}
Suppose that $(\gamma,\ell^+,\ell^-)$ and $(\overline{\gamma},\overline{\ell}^+,\overline{\ell}^-):I \to \R^3_1 \times \Delta_4$ are $(\bn,\overline{\bn})$-mates with $\overline{\gamma}(t)=\gamma(t)+\lambda(t)\bn(t)$. 
Then $\lambda$ is a non-zero constant. 
\end{proposition}
\demo
By differentiating $\overline{\gamma}(t)=\gamma(t)+\lambda(t)\bn(t)$, we have 
\begin{align*}
\dot{\overline{\gamma}}(t)
&=\overline{\alpha}(t) \overline{\ell}^+(t)+\overline{\beta}(t) \overline{\ell}^-(t)\\
&=(\alpha(t)+\lambda(t)\kappa_2(t))\ell^+(t)+(\beta(t)+\lambda(t)\kappa_3(t))\ell^-(t)+\dot{\lambda}(t)\bn(t).
\end{align*}
Since $\bn(t)=\overline{\bn}(t)$, we have $\dot{\lambda}(t)=0$ for all $t \in I$.
Hence $\lambda$ is a constant. 
If $\lambda$ is zero, then $\gamma$ and $\overline{\gamma}$ are the same. 
It follows that $\lambda$ is a non-zero constant. 
\enD
\begin{theorem}\label{spacelike-mate-1}
$(\gamma,\ell^+,\ell^-):I \to \R^3_1 \times \Delta_4$ is always a $(\bn,\overline{\bn})$-Bertrand lightcone framed curve. 
\end{theorem}
\demo
If we consider $(\overline{\gamma},\overline{\ell}^+,\overline{\ell}^-):I \to \R^3_1 \times \Delta_4$ by 
$$
\overline{\gamma}(t)=\gamma(t)+\lambda \bn(t), \ \overline{\ell}^+(t)=\ell^+(t), \ \overline{\ell}^-(t)=\ell^-(t),
$$
where $\lambda$ is a non-zero constant. 
By Proposition \ref{lambda}, 
$$
\dot{\overline{\gamma}}(t)=(\alpha(t)+\lambda(t)\kappa_2(t))\ell^+(t)+(\beta(t)+\lambda(t)\kappa_3(t))\ell^-(t).
$$
It follows that $(\overline{\gamma},\overline{\ell}^+,\overline{\ell}^-):I \to \R^3_1 \times \Delta_4$ is a lightcone framed curve with the curvature $(\kappa_1,\kappa_2,\kappa_3,\alpha+\lambda\kappa_2,\beta+\lambda\kappa_3)$.
Hence, $(\gamma,\ell^+,\ell^-)$ and $(\overline{\gamma},\overline{\ell}^+,\overline{\ell}^-)$ are $(\bn,\overline{\bn})$-mates.
\enD
\begin{remark}
{\rm 
If $(\gamma,\ell^+,\ell^-):I \to \R^3_1 \times \Delta_4$ is a $(\bn,\overline{\bn})$-Bertrand lightcone framed curve, then $\overline{\gamma}=\gamma+\lambda \bn$, where $\lambda$ is non-zero constant. 
Hence $\overline{\gamma}$ is a parallel curve of $\gamma$ with respect to $\bn$. 
We can also take $(\overline{\ell}^+,\overline{\ell}^-)$ as $(c \ell^+, (1/c)\ell^-)$ in Theorem \ref{spacelike-mate-1}, where $c:I \to \R$ is a non-zero function. See the special frame in \S 2.
}
\end{remark}

\begin{theorem}\label{spacelike-mate2}
$(\gamma,\ell^+,\ell^-):I \to \R^3_1 \times \Delta_4$ is a $(\bn^S,\overline{\bn}^S)$-Bertrand lightcone framed curve if and only if there exist smooth functions $\lambda,\theta:I \to \R$ with $\lambda \not\equiv 0$ such that 
$$
\lambda(t) (\kappa_2(t)-\kappa_3(t))\cosh \theta(t)+(\alpha(t)+\beta(t)+\lambda(t)\kappa_1(t))\sinh \theta(t)=0 
$$
for all $t \in I$. 
\end{theorem}
\demo
Suppose that $(\gamma,\ell^+,\ell^-):I \to \R^3_1 \times \Delta_4$ is a $(\bn^S,\overline{\bn}^S)$-Bertrand lightcone framed curve. Then there exists another lightcone framed curve $(\overline{\gamma},\overline{\ell}^+,\overline{\ell}^-)$ such that $(\gamma,\ell^+,\ell^-)$ and $(\overline{\gamma},\overline{\ell}^+,\overline{\ell}^-)$ are $(\bn^S,\overline{\bn}^S)$-mates, that is, 
there exists a smooth function $\lambda:I \to \R$ with $\lambda \not\equiv 0$ such that $\overline{\gamma}(t)=\gamma(t)+\lambda(t)\bn^S(t)$ and $\bn^S(t)=\overline{\bn}^S(t)$ for all $t \in I$. 
By differentiating $\overline{\gamma}(t)=\gamma(t)+\lambda(t) \bn^S(t)$, we have
$$
\overline{a}(t)\overline{\bn}^T(t)+\overline{b}(t)\overline{\bn}^S(t)=(a(t)+\lambda(t)\kappa_1(t))\bn^T(t)+(b(t)+\dot{\lambda}(t))\bn^S(t)-\lambda(t)\kappa^S(t)\bn(t).
$$
Since $\bn^S(t)=\overline{\bn}^S(t)$, we have $\overline{b}(t)=b(t)+\dot{\lambda}(t)$. 
If  we denote 
$$
\overline{\bn}(t)=a_1(t)\bn(t)+b_1(t)\bn^T(t), \ \overline{\bn}^T(t)=a_2(t)\bn(t)+b_2(t)\bn^T(t),
$$
where $a_1,a_2,b_1b_2:I \to \R$ are smooth functions, then $a_1^2(t)-b_1^2(t)=1, a_2^2(t)-b_2^2(t)=-1$ and $a_1(t)a_2(t)-b_1(t)b_2(t)=0$. 
Since $\overline{\bn}^S(t)=\overline{\bn}(t) \wedge \overline{\bn}^T(t)=(a_1(t)b_2(t)-a_2(t)b_1(t))\bn^S(t)$, $a_1(t)b_2(t)-a_2(t)b_1(t)=1$. 
It follows that $a_2(t)=b_1(t), b_2(t)=a_1(t)$.
By  $a_1^2(t)-b_1^2(t)=1$, there exists a smooth function $\theta:I \to \R$ such that 
\begin{align*}
\begin{pmatrix}
\overline{\bn}(t)\\
\overline{\bn}^T(t)
\end{pmatrix}
= \pm
\begin{pmatrix}
\cosh \theta(t) & \sinh \theta(t)\\
\sinh \theta(t) & \cosh \theta(t)
\end{pmatrix}
\begin{pmatrix}
{\bn}(t)\\
{\bn}^T(t)
\end{pmatrix}.
\end{align*}
Then we have $\pm \overline{a}(t)\sinh \theta(t)=-\lambda(t)\kappa^S(t)$ and 
$\pm \overline{a}(t)\cosh \theta(t)=a(t)+\lambda(t)\kappa_1(t)$.
It follows that $\overline{a}(t)=\pm (\lambda(t)\kappa^S(t) \sinh  \theta(t)+(a(t)+\lambda(t)\kappa_1(t))\cosh \theta(t))$ and $\lambda(t)\kappa^S(t) \cosh \theta(t)+(a(t)+\lambda(t)\kappa_1(t))\sinh \theta(t)=0$.
Therefore, we have
$$
\lambda(t) (\kappa_2(t)-\kappa_3(t))\cosh \theta(t)+(\alpha(t)+\beta(t)+\lambda(t)\kappa_1(t))\sinh  \theta(t)=0
$$
for all $t \in I$.
\par
Conversely, suppose that there exist smooth functions $\lambda,\theta:I \to \R$ with $\lambda \not\equiv 0$ such that $\lambda(t)\kappa^S(t) \cosh \theta(t)+(a(t)+\lambda(t)\kappa_1(t))\sinh \theta(t)=0$ for all $t \in I$. 
If we consider $\overline{\gamma}(t)=\gamma(t)+\lambda(t)\bn^S(t), \overline{\bn}(t)=\pm(\cosh  \theta(t)\bn(t)+\sinh \theta(t)\bn^T(t)), \overline{\bn}^T(t)=\pm (\sinh \theta(t)\bn(t)+\cosh \theta(t)\bn^T(t))$, 
then $\overline{\bn}^S(t)=\bn^S(t)$ and  
\begin{align*}
\dot{\overline{\gamma}}(t) &=(a(t)+\lambda(t)\kappa_1(t))\bn^T(t)+(b(t)+\dot{\lambda}(t))\bn^S(t)-\lambda(t)\kappa^S(t)\bn(t)\\
&=\pm (a(t)+\lambda(t)\kappa_1(t))(-\sinh \theta(t) \overline{\bn}(t)+\cosh \theta(t)\overline{\bn}^T(t))+(b(t)+\dot{\lambda}(t))\overline{\bn}^S(t)\\
&\quad \mp \lambda(t)\kappa^S(t)(\cosh \theta(t) \overline{\bn}(t)-\sinh \theta(t) \overline{\bn}^T(t)) \\
&=\pm (\lambda(t)\kappa^S(t) \sinh  \theta(t)+(a(t)+\lambda(t)\kappa_1(t))\cosh \theta(t))\overline{\bn}^T(t)\\
&\quad+(b(t)+\dot{\lambda}(t))\overline{\bn}^S(t)\\
&= \overline{a}(t)\overline{\bn}^T(t)+\overline{b}(t)\overline{\bn}^S(t)\\
&=\frac{1}{2}(\overline{a}(t)+\overline{b}(t))\overline{\ell}^+(t)+\frac{1}{2}(\overline{a}(t)-\overline{b}(t))\overline{\ell}^-(t),
\end{align*}
where 
\begin{align*}
& \overline{a}(t)=\pm (\lambda(t)\kappa^S(t) \sinh  \theta(t)+(a(t)+\lambda(t)\kappa_1(t))\cosh \theta(t)), \ \overline{b}(t)=b(t)+\dot{\lambda}(t),\\
& \overline{\ell}^+(t)=\overline{\bn}^T(t)+\overline{\bn}^S(t), \ \overline{\ell}^-(t)=\overline{\bn}^T(t)-\overline{\bn}^S(t).
\end{align*}
It follows that $(\overline{\gamma},\overline{\ell}^+,\overline{\ell}^-):I \to \R^3_1 \times \Delta_4$ is a lightcone framed curve. 
Hence, $(\gamma,\ell^+,\ell^-)$ and $(\overline{\gamma},\overline{\ell}^+,\overline{\ell}^-)$ are $(\bn^S,\overline{\bn}^S)$-mates.
\enD

\begin{proposition}\label{curvature-spacelike-mate2}
Suppose that there exist smooth functions $\lambda,\theta:I \to \R$ with $\lambda \not\equiv 0$ such that
$$
\lambda(t) (\kappa_2(t)-\kappa_3(t))\cosh \theta(t)+(\alpha(t)+\beta(t)+\lambda(t)\kappa_1(t))\sinh \theta(t)=0 
$$
for all $t \in I$. 
Then the curvature of the lightcone framed curve $(\overline{\gamma},\overline{\ell}^+,\overline{\ell}^-):I \to \R^3_1 \times \Delta_4$, $\overline{\gamma}(t)=\gamma(t)+\lambda(t)\bn^S(t),\overline{\ell}^+(t)=\overline{\bn}^T(t)+\overline{\bn}^S(t), \overline{\ell}^-(t)=\overline{\bn}^T(t)-\overline{\bn}^S(t), 
\overline{\bn}^{T}(t)=\pm(\sinh \theta(t) \bn(t)+\cosh \theta(t)\bn^T(t)), 
\overline{\bn}^S(t)=\bn^S(t)$ is given by 
$$
(\overline{\kappa}_1,\overline{\kappa}_2,\overline{\kappa}_3,\overline{\alpha},\overline{\beta})=(\overline{\kappa}_1,(1/2)(\overline{\kappa}^T+\overline{\kappa}^S),(1/2)(\overline{\kappa}^T-\overline{\kappa}^S),(1/2)(\overline{a}+\overline{b}),(1/2)(\overline{a}-\overline{b})),
$$ 
where
\begin{align*}
\overline{\kappa}_1(t)&=\pm(\kappa_1(t)\cosh \theta(t)+\kappa^S(t) \sinh \theta(t)), \\ 
\overline{\kappa}^T(t)&=\pm(\dot{\theta}(t)+\kappa^T(t)), \\ 
\overline{\kappa}^S(t)&=\pm(\kappa_1(t)\sinh \theta(t)+\kappa^S(t)\cosh \theta(t)),\\
\overline{a}(t)&=\pm(\lambda(t)\kappa^S(t) \sinh \theta(t)+(a(t)+\lambda(t)\kappa_1(t))\cosh \theta(t)), \\ 
\overline{b}(t)&=b(t)+\dot{\lambda}(t).
\end{align*}
\end{proposition}
\demo
By a direct calculation,  
\begin{align*}
\dot{\overline{\bn}^T}(t)&=\pm (\dot{\theta}(t)\cosh \theta(t) \bn(t)+\sinh \theta(t)(\kappa^T(t)\bn^T(t)+\kappa^S(t)\bn^S(t))\\
&\quad +\dot{\theta}(t)\sinh \theta(t)\bn^T(t)+\cosh \theta(t)(\kappa_1(t)\bn^T(t)-\kappa^S(t)\bn(t))) \\
&=\pm (\kappa_1 (t)\cosh \theta(t)+\kappa^S(t) \sinh \theta(t))\overline{\bn}^{S}\pm(\dot{\theta}(t)+\kappa^T(t))\overline{\bn}(t), \\
\dot{\overline{\bn}^S}(t)&=\kappa_1(t)\bn^S(t)-\kappa^S(t)\bn(t)\\
&=\pm(\kappa_1 (t)\cosh \theta(t)+\kappa^S(t)\sinh(t))\overline{\bn}^T(t)\\
&\quad\mp(\kappa_1(t)\sinh \theta(t)+\kappa^S(t)\sinh \theta(t))\overline{\bn}(t).
\end{align*}
Then we have 
\begin{align*}
&\overline{\kappa}_1(t)=\pm(\kappa_1(t) \cosh \theta(t)+\kappa^S(t)\sinh(t)), \ 
\overline{\kappa}^T(t)=\pm(\dot{\theta}(t)+\kappa^T(t)), \\ 
&\overline{\kappa}^S(t)=\pm(\kappa_1(t)\sinh \theta(t)+\kappa^S(t)\sinh \theta(t)).
\end{align*}
By the proof of Theorem \ref{spacelike-mate2}, we also have 
$$
\overline{a}(t)=\pm (\lambda(t)\kappa^S(t) \sinh  \theta(t)+(a(t)+\lambda(t)\kappa_1(t))\cosh \theta(t)), \ 
\overline{b}(t)=b(t)+\dot{\lambda}(t).
$$
\enD


\begin{theorem}\label{spacelike-mate3}
$(\gamma,\ell^+,\ell^-):I \to \R^3_1 \times \Delta_4$ is a $(\bn,\overline{\bn}^S)$-Bertrand lightcone framed curve if and only if there exist smooth functions $\lambda,\theta:I \to \R$ with $\lambda \not\equiv 0$ such that 
\begin{align*}
&\left(\alpha(t)+\beta(t)+\lambda(t)(\kappa_2(t)+\kappa_3(t))\right)\sinh\theta(t)\\
&\quad-\left(\alpha(t)-\beta(t)+\lambda(t)(\kappa_2(t)-\kappa_3(t))\right)\cosh\theta(t)=0 
\end{align*}
for all $t \in I$. 
\end{theorem}
\demo
Suppose that $(\gamma,\ell^+,\ell^-):I \to \R^3_1 \times \Delta_4$ is a $(\bn,\overline{\bn}^S)$-Bertrand lightcone framed curve. Then there exists another lightcone framed curve $(\overline{\gamma},\overline{\ell}^+,\overline{\ell}^-)$ such that $(\gamma,\ell^+,\ell^-)$ and $(\overline{\gamma},\overline{\ell}^+,\overline{\ell}^-)$ are $(\bn,\overline{\bn}^S)$-mates, that is, 
there exists a smooth function $\lambda:I \to \R$ with $\lambda \not\equiv 0$ such that $\overline{\gamma}(t)=\gamma(t)+\lambda(t)\bn(t)$ and $\bn(t)=\overline{\bn}^S(t)$ for all $t \in I$. 
By differentiating $\overline{\gamma}(t)=\gamma(t)+\lambda(t) \bn(t)$, we have
$$
\overline{a}(t)\overline{\bn}^T(t)+\overline{b}(t)\overline{\bn}^S(t)=(a(t)+\lambda(t)\kappa^T(t))\bn^T(t)+(b(t)+\lambda(t)\kappa^S(t)
)\bn^S(t)+\dot{\lambda}(t)\bn(t).
$$
Since $\bn(t)=\overline{\bn}^S(t)$, we have $\overline{b}(t)=\dot{\lambda}(t)$. 
If we denote 
$$
\overline{\bn}(t)=a_1(t)\bn^T(t)+b_1(t)\bn^S(t), \ \overline{\bn}^T(t)=a_2(t)\bn^T(t)+b_2(t)\bn^S(t),
$$
where $a_1,a_2,b_1b_2:I \to \R$ are smooth functions, then $-a_1^2(t)+b_1^2(t)=1, -a_2^2(t)+b_2^2(t)=-1$ and $-a_1(t)a_2(t)+b_1(t)b_2(t)=0$. 
Since $\overline{\bn}^S(t)=\overline{\bn}(t) \wedge \overline{\bn}^T(t)=(a_1(t)b_2(t)-a_2(t)b_1(t))\bn(t)$, $a_1(t)b_2(t)-a_2(t)b_1(t)=1$. 
It follows that $a_2(t)=b_1(t), b_2(t)=a_1(t)$.
By $-a_1^2(t)+b_1^2(t)=1$, there exists a smooth function $\theta:I \to \R$ such that 
\begin{align*}
\begin{pmatrix}
\overline{\bn}(t)\\
\overline{\bn}^T(t)
\end{pmatrix}
= \pm
\begin{pmatrix}
\sinh \theta(t) & \cosh \theta(t)\\
-\cosh \theta(t) & -\sinh \theta(t)
\end{pmatrix}
\begin{pmatrix}
{\bn}^T(t)\\
{\bn}^S(t)
\end{pmatrix}.
\end{align*}
Then we have $\mp \overline{a}(t)\cosh \theta(t)=a(t)+\lambda(t)\kappa^T(t)$ and 
$\mp \overline{a}(t)\sinh \theta(t)=b(t)+\lambda(t)\kappa^S(t)$.
It follows that $\overline{a}(t)=\mp ((a+\lambda(t)\kappa^T(t)) \cosh \theta(t)-(b(t)+\lambda(t)\kappa^S(t))\\\sinh \theta(t))$ and $(a+\lambda(t)\kappa^T(t)) \sinh \theta(t)-(b(t)+\lambda(t)\kappa^S(t))\cosh \theta(t)=0$.
Therefore, we have
\begin{align*}
&(\alpha(t)+\beta(t)+\lambda(t)(\kappa_2(t)+\kappa_3(t)))\sinh \theta(t)\\
&\quad-(\alpha(t)-\beta(t)+\lambda(t)(\kappa_2(t)-\kappa_3(t)))\cosh \theta(t)=0
\end{align*}
for all $t \in I$.
\par
Conversely, suppose that there exist smooth functions $\lambda,\theta:I \to \R$ with $\lambda \not\equiv 0$ such that $(a+\lambda(t)\kappa^T(t)) \sinh \theta(t)-(b(t)+\lambda(t)\kappa^S(t))\cosh \theta(t)=0$ for all $t \in I$. 
If we consider $\overline{\gamma}(t)=\gamma(t)+\lambda(t)\bn(t), \overline{\bn}^T(t)=\mp(\cosh \theta(t)\bn^T(t)+\sinh \theta(t)\bn^S(t)), \overline{\bn}(t)=\pm (\sinh \theta(t)\bn^T(t)+\cosh \theta(t)\bn^S(t))$, 
then $\overline{\bn}^S(t)=\bn(t)$ and  
\begin{align*}
\dot{\overline{\gamma}}(t) &=(a(t)+\lambda(t)\kappa^T(t))\bn^T(t)+(b(t)+\lambda(t)\kappa^S(t))\bn^S(t)+\dot{\lambda}(t)\bn(t)\\
&=\mp (a(t)+\lambda(t)\kappa^T(t))(\sinh \theta(t) \overline{\bn}(t)+\cosh \theta(t)\overline{\bn}^T(t))+\dot{\lambda}(t)\overline{\bn}^S(t)\\
&\quad \pm (b(t)+\lambda(t)\kappa^S(t))(\cosh \theta(t) \overline{\bn}(t)+\sinh \theta(t) \overline{\bn}^T(t)) \\
&=\mp ((a+\lambda(t)\kappa^T(t)) \cosh \theta(t)-(b(t)+\lambda(t)\kappa^S(t))\sinh \theta(t))\overline{\bn}^T(t)\\
&\quad+\dot{\lambda}(t)\overline{\bn}^S(t) \\
&=\overline{a}(t)\overline{\bn}^T(t)+\overline{b}(t)\overline{\bn}^S(t)\\
&=\frac{1}{2}(\overline{a}(t)+\overline{b}(t))\overline{\ell}^+(t)+\frac{1}{2}(\overline{a}(t)-\overline{b}(t))\overline{\ell}^-(t),
\end{align*}
where 
\begin{align*}
& \overline{a}(t)=\mp ((a+\lambda(t)\kappa^T(t)) \cosh \theta(t)-(b(t)+\lambda(t)\kappa^S(t))\sinh \theta(t)), \ \overline{b}(t)=\dot{\lambda}(t),\\
& \overline{\ell}^+(t)=\overline{\bn}^T(t)+\overline{\bn}^S(t), \ \overline{\ell}^-(t)=\overline{\bn}^T(t)-\overline{\bn}^S(t).
\end{align*}
It follows that $(\overline{\gamma},\overline{\ell}^+,\overline{\ell}^-):I \to \R^3_1 \times \Delta_4$ is a lightcone framed curve. 
Hence, $(\gamma,\ell^+,\ell^-)$ and $(\overline{\gamma},\overline{\ell}^+,\overline{\ell}^-)$ are $(\bn,\overline{\bn}^S)$-mates.
\enD

\begin{proposition}\label{curvature-spacelike-mate3}
Suppose that there exist smooth functions $\lambda,\theta:I \to \R$ with $\lambda \not\equiv 0$ such that
\begin{align*}
&\left(\alpha(t)+\beta(t)+\lambda(t)(\kappa_2(t)+\kappa_3(t))\right)\sinh\theta(t)\\
&\quad-\left(\alpha(t)-\beta(t)+\lambda(t)(\kappa_2(t)-\kappa_3(t))\right)\cosh\theta(t)=0 
\end{align*}
for all $t \in I$. 
Then the curvature of the lightcone framed curve $(\overline{\gamma},\overline{\ell}^+,\overline{\ell}^-):I \to \R^3_1 \times \Delta_4$, $\overline{\gamma}(t)=\gamma(t)+\lambda(t)\bn(t),\overline{\ell}^+(t)=\overline{\bn}^T(t)+\overline{\bn}^S(t), \overline{\ell}^-(t)=\overline{\bn}^T(t)-\overline{\bn}^S(t), 
\overline{\bn}^{T}(t)=\mp(\cosh \theta(t) \bn^T(t)+\sinh \theta(t)$$\bn^S(t)), \overline{\bn}^S(t)=\bn(t)$ is given by 
$$
(\overline{\kappa}_1,\overline{\kappa}_2,\overline{\kappa}_3,\overline{\alpha},\overline{\beta})=(\overline{\kappa}_1,(1/2)(\overline{\kappa}^T+\overline{\kappa}^S),(1/2)(\overline{\kappa}^T-\overline{\kappa}^S),(1/2)(\overline{a}+\overline{b}),(1/2)(\overline{a}-\overline{b})),
$$ 
where
\begin{align*}
\overline{\kappa}_1(t)&=\mp\left(\kappa^T(t)\cosh \theta(t)-\kappa^S(t) \sinh \theta(t)\right), \\ 
\overline{\kappa}^T(t)&=\mp(\dot\theta(t)+\kappa_1(t)), \\ 
\overline{\kappa}^S(t)&=\pm\left(\kappa^T(t)\sinh \theta(t)-\kappa^S(t)\cosh \theta(t)\right),\\
\overline{a}(t)&=\mp\left((a(t)+\lambda(t)\kappa^T(t) )\cosh \theta(t)-(b(t)+\lambda(t)\kappa^S(t))\sinh \theta(t)\right), \\ 
\overline{b}(t)&=\dot{\lambda}(t).
\end{align*}
\end{proposition}
\demo
By a direct calculation,  
\begin{align*}
\dot{\overline{\bn}^T}(t)&=\mp (\dot{\theta}(t)\sinh \theta(t) \bn^T(t)+\cosh \theta(t)(\kappa_1(t)\bn^S(t)+\kappa^T(t)\bn(t))\\
&\quad +\dot{\theta}(t)\cosh \theta(t)\bn^S(t)+\sinh \theta(t)(\kappa_1(t)\bn^T(t)-\kappa^S(t)\bn(t))) \\
&=\mp (\kappa^T (t)\cosh \theta(t)-\kappa^S(t) \sinh \theta(t))\overline{\bn}^{S}(t)\mp(\dot{\theta}(t)+\kappa_1(t))\overline{\bn}(t), \\
\dot{\overline{\bn}}(t)&=\pm (\dot{\theta}(t)\cosh \theta(t) \bn^T(t)+\sinh \theta(t)(\kappa_1(t)\bn^S(t)+\kappa^T(t)\bn(t))\\
&\quad +\dot{\theta}(t)\sinh \theta(t)\bn^S(t)+\cosh \theta(t)(\kappa_1(t)\bn^T(t)-\kappa^S(t)\bn(t))) \\
&=\pm (\kappa^T (t)\sinh \theta(t)-\kappa^S(t) \cosh \theta(t))\overline{\bn}^{S}(t)\mp(\dot{\theta}(t)+\kappa_1(t))\overline{\bn}^{T}(t).
\end{align*}
Then we have 
\begin{align*}
&\overline{\kappa}_1(t)=\mp\left(\kappa^T(t)\cosh \theta(t)-\kappa^S(t) \sinh \theta(t)\right), \ 
\overline{\kappa}^T(t)=\mp(\dot\theta(t)+\kappa_1(t)), \\ 
&\overline{\kappa}^S(t)=\pm\left(\kappa^T(t)\sinh \theta(t)-\kappa^S(t)\cosh \theta(t)\right).
\end{align*}
By the proof of Theorem \ref{spacelike-mate2}, we also have 
$$
\overline{a}(t)=\mp ((a+\lambda(t)\kappa^T(t)) \cosh \theta(t)-(b(t)+\lambda(t)\kappa^S(t))\sinh \theta(t)), \
\overline{b}(t)=\dot{\lambda}(t).
$$
\enD

\begin{theorem}\label{spacelike-mate4}
$(\gamma,\ell^+,\ell^-):I \to \R^3_1 \times \Delta_4$ is always a $(\bn^S,\overline{\bn})$-Bertrand lightcone framed curve. 
\end{theorem}
\demo
If we consider $(\overline{\gamma},\overline{\ell}^+,\overline{\ell}^-):I \to \R^3_1 \times \Delta_4$ by 
\begin{align*}
&\overline{\gamma}(t)=\gamma(t)+\lambda(t) \bn^S(t), \ \overline{\ell}^+(t)=\widetilde{\ell}^-(t)=\bn^T(t)+\bn(t), \\ 
&\overline{\ell}^-(t)=\widetilde{\ell}^+(t)=\bn^T(t)-\bn(t),
\end{align*}
where $\lambda(t)=-\int b(t) dt$. By differentiating $\overline{\gamma}(t)=\gamma(t)-(\int b(t) dt) \bn^S(t)$, we have
then we have
$$
\dot{\overline{\gamma}}(t)=\left(a(t)-\kappa_1(t)\int b(t) dt\right)\overline{\bn}^T(t)+\left(-\kappa^S(t)\int b(t) dt\right)\overline{\bn}^S(t).
$$
It follows that $(\overline{\gamma},\overline{\ell}^+,\overline{\ell}^-):I \to \R^3_1 \times \Delta_4$ is a lightcone framed curve with the curvature 
$$
(\overline{\kappa}_1,\overline{\kappa}_2,\overline{\kappa}_3,\overline{\alpha},\overline{\beta})=(\overline{\kappa}_1,(1/2)(\overline{\kappa}^T+\overline{\kappa}^S),(1/2)(\overline{\kappa}^T-\overline{\kappa}^S),(1/2)(\overline{a}+\overline{b}),(1/2)(\overline{a}-\overline{b})),
$$ 
where
\begin{align*}
&\overline{\kappa}_1(t)=-\kappa^T(t), \ \overline{\kappa}^T(t)=\kappa^S(t), \ \overline{\kappa}^S(t)=\kappa_1(t), \\
&\overline{a}(t)=a(t)-\kappa_1(t)\int b(t) dt, \ \overline{b}(t)=-\kappa^S(t)\int b(t) dt.  
\end{align*}
Hence, $(\gamma,\ell^+,\ell^-)$ and $(\overline{\gamma},\overline{\ell}^+,\overline{\ell}^-)$ are $(\bn^S,\overline{\bn})$-mates.
\enD
\begin{remark}
{\rm 
If $(\gamma,\ell^+,\ell^-):I \to \R^3_1 \times \Delta_4$ is a $(\bn^S,\overline{\bn})$-Bertrand lightcone framed curve, then $\overline{\gamma}(t)=\gamma(t)-(\int b(t) dt) \bn^S(t)$. 
Hence, $\overline{\gamma}$ is a psedo-circular involute of $\gamma$ (cf. \cite{Pei-Takahashi-Zhang, Zhang-Li-Pei}). 
}
\end{remark}

\section{Examples}
We give concrete example of Bertrand lightcone framed curves.
\begin{example}
{\rm 
Let $p, q\in \R$ and $n, m\in \N$. Suppose that $n\neq m$ (respectively, $n= m$).
Then we define $(\gamma,\ell^+,\ell^-):I \to \R^3_1 \times \Delta_4$ by
\begin{gather*}
\gamma(t)=\Bigg(-\frac{p}{m}\cos mt, -\frac{q}{2}\left(\frac{\cos (m-n)t}{m-n}+\frac{\cos(m+n)t}{m+n}\right),\\
-\frac{q}{2}\left(\frac{\sin (m-n)t}{m-n}-\frac{\sin(m+n)t}{m+n}\right)\Bigg) \\
\left({\rm respectively}, \gamma(t)=\left(-\frac{p}{m}\cos mt, -\frac{q}{4m}\cos 2mt, \frac{q}{2}t-\frac{q}{4m}\sin 2mt\right)\right),\\
\ell^+(t)=(1, \cos nt, \sin nt),\\
\ell^-(t)=(1, -\cos nt, -\sin nt).
\end{gather*}
Then $(\gamma,\ell^+,\ell^-):I \to \R^3_1 \times \Delta_4$ is a lightcone framed curve with the curvature 
\begin{align*}
{\kappa}_1(t)=0, \ {\kappa}_2(t)=-\frac{n}{2}, \ {\kappa}_3=\frac{n}{2}, \ {\alpha}(t)=\frac{p+q}{2}\sin mt, \ 
{\beta}(t)=\frac{p-q}{2}\sin mt.
\end{align*}
By definition, we have
$$
{\kappa}^T(t)=0, \ {\kappa}^S(t)=-n, \ {a}(t)=p\sin mt, \ 
{b}(t)=q\sin mt.
$$
If we take $\lambda(t)=(k(t)/n)(p-q)\sin mt$ (respectively, $\lambda(t)=-(k(t)/n)(p+q)\sin mt$) with $\lambda \not\equiv 0$, then the condition of Theorem \ref{lightlike-mate-1} (respectively, of Theorem \ref{lightlike-mate-3}) is satisfied. 
Hence, $(\gamma,\ell^+,\ell^-)$ is a $(\ell^+, \overline{\ell}^+)$ and $(\ell^+, \overline{\ell}^-)$ (respectively, $(\ell^-, \overline{\ell}^+)$ and $(\ell^-, \overline{\ell}^-)$)-Bertrand lightcone framed curve. 
Similarly, if we take $\lambda(t)=((k(t)p+q)/n)\sin mt$ (respectively, $\lambda(t)=-((k(t)p-q)/n)\sin mt$) with $\lambda \not\equiv 0$, then the condition of Theorem \ref{lightlike-mate-5} (respectively, of Theorem \ref{lightlike-mate-7}) is satisfied. 
Hence, $(\gamma,\ell^+,\ell^-)$ is a $(\widetilde{\ell}^+,\overline{\ell}^+)$ and $(\widetilde{\ell}^+,\overline{\ell}^-)$ (respectively, $(\widetilde{\ell}^-, \overline{\ell}^+)$ and $(\widetilde{\ell}^-, \overline{\ell}^-)$)-Bertrand lightcone framed curve. \par
If we take $\theta(t)=n\pi$, then the condition of Theorem \ref{timelike-mate} is satisfied.
Hence, $(\gamma,\ell^+,\ell^-)$ is a $(\bn^T,\overline{\bn}^T)$-Bertrand lightcone framed curve. \par
If we take $\lambda(t)=(p\tanh \theta(t)/n)\sin mt$ (respectively, $\lambda(t)=((p\tanh \theta(t)-q)/n)\sin mt$) with $\lambda \not\equiv 0$, then the condition of Theorem \ref{spacelike-mate2} (respectively, of Theorem \ref{spacelike-mate3}) is satisfied. 
Hence, $(\gamma,\ell^+,\ell^-)$ is a $(\bn^S,\overline{\bn}^S)$ (respectively, $(\bn,\overline{\bn}^S)$)-Bertrand lightcone framed curve. 
Moreover, if we take $\lambda(t)=(q/m)\\\cos mt \not\equiv 0$, then the condition of Theorem \ref{spacelike-mate4} is satisfied. 
Hence, $(\gamma,\ell^+,\ell^-)$ is a $(\bn^S,\overline{\bn})$-Bertrand lightcone framed curve. 
Similarly, if $\lambda$ is a non-zero constant, then the condition of Theorem \ref{spacelike-mate-1} is satisfied. 
Hence, $(\gamma,\ell^+,\ell^-)$ is a $(\bn,\overline{\bn})$-Bertrand lightcone framed curve. 
}
\end{example}
\begin{example}
{\rm 
Let $(\gamma,\ell^+,\ell^-):I \to \R^3_1 \times \Delta_4$ be
\begin{align*}
&\gamma(t)=\Biggl(\int \alpha(t)\left(\int \kappa_3(t)dt\right)^2 dt+\int (\alpha(t)+\beta(t)) dt, \\
&\quad-\int \alpha(t)\left(\int \kappa_3(t)dt\right)^2 dt+\int (\alpha(t)-\beta(t)) dt,
2\left(\int\alpha(t)\left(\int\kappa_3(t)\right)dt\right)dt\Biggr),\\
&\ell^+(t)=\Biggl(\left(\int \kappa_3(t)dt\right)^2+1, \ -\left(\int \kappa_3(t)dt\right)^2+1, \ 2\int \kappa_3(t)dt\Biggr),\\
&\ell^-(t)=(1, \ -1, \ 0).
\end{align*}
Then $(\gamma,\ell^+,\ell^-):I \to \R^3_1 \times \Delta_4$ is a lightcone framed curve with the curvature $(0, 0, {\kappa}_3, {\alpha}, {\beta}).$
By definition, we have ${\kappa}^T(t)={\kappa}_3(t), {\kappa}^S(t)=-{\kappa}_3,  {a}(t)=\alpha(t)+\beta(t), \ {b}(t)=\alpha(t)-\beta(t)$. 
Since $\kappa_2(t)=0$, then the condition of Theorem \ref{lightlike-mate-3} is satisfied. Hence, $(\gamma,\ell^+,\ell^-)$ is always a $(\ell^-,\overline{\ell}^+)$-Bertrand lightcone framed curve. \par
If we take ${\kappa}_3(t)=\sin t$, $\alpha(t)=\sin t\cos t$ and $\beta(t)=\cos t$, then $(\gamma,\ell^+,\ell^-)$ is given by 
\begin{align*}
\gamma(t)&=\left(-\frac{1}{4}\cos^4 t-\frac{1}{2}\cos^2 t+\sin t, \ \frac{1}{4}\cos^4 t-\frac{1}{2}\cos^2 t-\sin t, \ \frac{2}{3}\cos^3 t\right),\\
\ell^+(t)&=\left(\cos^2 t+1, \ -\cos^2 t+1, \ -2\cos t \right),\\
\ell^-(t)&=(1, \ -1, \ 0),
\end{align*}
up to constant. 
Since $\alpha(t)\beta(t)=\sin t\cos^2 t$, $\gamma$ is a mixed type curve with singular points at $t=\pi/2$ and $3\pi/2$.
We also have $\kappa^T(t)=\sin t$, $\kappa^S(t)=-\sin t$, $a(t)=\cos t(\sin t+1)$ and $b(t)=\cos t(\sin t-1)$. 
Since $\beta(t)=\cos t\not\equiv0$, $(\gamma,\ell^+,\ell^-)$ is always a $(\ell^+,\overline{\ell}^+)$-Bertrand lightcone framed curve by Remark \ref{lightlike-mate-1-re} $(2)$. 
If we take $k(t)=1$ and $\lambda(t)=-2\cos t$ (respectively, $k(t)=-1$ and $\lambda(t)=-2\cos t$), then the condition of Theorem \ref{lightlike-mate-5} (respectively, of Theorem \ref{lightlike-mate-7}) is satisfied. Hence, $(\gamma,\ell^+,\ell^-)$ is a $(\widetilde{\ell}^+,\overline{\ell}^+)$ and $(\widetilde{\ell}^+,\overline{\ell}^-)$ (respectively, $(\widetilde{\ell}^-,\overline{\ell}^+)$ and $(\widetilde{\ell}^-,\overline{\ell}^-)$)-Bertrand lightcone framed curve. 
Similarly, if we take $\theta(t)=t$ and $\lambda(t)=\sin t-1$, then the condition of Theorem \ref{timelike-mate} is satisfied. Hence, $(\gamma,\ell^+,\ell^-)$ is a $(\bn^T,\overline{\bn}^T)$-Bertrand lightcone framed curve. \par
If we take ${\kappa}_3(t)=\sinh t$, $\alpha(t)=\sinh t\cosh t$ and $\beta(t)=-\sinh t$, then $(\gamma,\ell^+,\ell^-)$ is given by 
\begin{align*}
&\gamma(t)=\left(\frac{\cosh^4 t}{4}+\frac{\cosh^2 t}{2}-\cosh t, \ -\frac{\cosh^4 t}{4}+\frac{\cosh^2 t}{2}+\cosh t, \ \frac{2}{3}\cosh^3 t\right),\\
&\ell^+(t)=\left(\cosh^2 t+1, \ -\cosh^2 t+1, \ 2\cosh t \right),\\
&\ell^-(t)=(1, \ -1, \ 0),
\end{align*}
up to constant. 
Since $\alpha(t)\beta(t)=-\sinh^2 t\cosh t$, $\gamma$ is a spcelike curve with a singular point at $t=0$. 
We also have $\kappa^T(t)=\sinh t$, $\kappa^S(t)=-\sinh t$, $a(t)=\sinh t(\cosh t-1)$ and $b(t)=\sinh t(\cosh t+1)$. If we take a smooth function $\theta : I\to \R$ such that $\lambda(t)=e^{-2\theta(t)}\cosh t+1$, then the condition of Theorem \ref{spacelike-mate3} is satisfied. Hence, $(\gamma,\ell^+,\ell^-)$ is a $(\bn,\overline{\bn}^S)$-Bertrand lightcone framed curve. Similarly, if we take a smooth function $\theta : I\to \R$ such that $\lambda(t)=-\tanh \theta(t)(\cosh t-1)$ and $\theta\not\equiv0$, then the condition of Theorem \ref{spacelike-mate2} is satisfied. Hence, $(\gamma,\ell^+,\ell^-)$ is a $(\bn^S,\overline{\bn}^S)$-Bertrand lightcone framed curve.
}
\end{example}

Nozomi Nakatsuyama, 
\\
Muroran Institute of Technology, Muroran 050-8585, Japan,
\\
E-mail address: 23043042@muroran-it.ac.jp
\\
\\
Masatomo Takahashi, 
\\
Muroran Institute of Technology, Muroran 050-8585, Japan,
\\
E-mail address: masatomo@muroran-it.ac.jp

\end{document}